\pdfoutput=1
\documentclass[11pt,a4paper]{scrartcl}  
\usepackage[l2tabu, orthodox]{nag}  
\usepackage[utf8]{inputenc}
\usepackage[pdfusetitle]{hyperref}

\usepackage{amsmath, amsthm, amssymb}

\usepackage{mathptmx}
\usepackage[scaled=.90]{helvet}
\usepackage{courier}

\usepackage{graphicx}
\usepackage{authblk}
\usepackage[semicolon, round]{natbib}
\usepackage{mathtools}  
\usepackage[usenames,dvipsnames]{color}  
\usepackage{listings}
\usepackage{tabu}

\usepackage{stmaryrd}  

\lstdefinelanguage[firedrake]{python}[]{python}{%
  emph={[2]Function,Mesh,FunctionSpace,interpolate,TrialFunction,TestFunction,Expression,DirichletBC,solve,assemble,assign,apply,FiniteElement,
  TensorProductElement,EnrichedElement,HDiv,HCurl,triangle,interval},
  emph={grad,dx,inner}
}

\definecolor{DarkBlue}{rgb}{0.00,0.00,0.55}
\definecolor{DarkRed}{rgb}{0.55,0.00,0.00}
\definecolor{DarkGreen}{rgb}{0.00,0.55,0.00}
\numberwithin{equation}{section}

\newcommand{\dx}{\,\mathrm{d} x}
\newcommand{\ints}[1]{\int_{#1}}  
\renewcommand{\P}[1]{{\mathrm{P}}_{#1}}
\newcommand{\DP}[1]{{\mathrm{DP}}_{#1}}
\newcommand{\Q}[1]{{\mathrm{Q}}_{#1}}
\newcommand{\DQ}[1]{{\mathrm{DQ}}_{#1}}
\newcommand{\B}[1]{{\mathrm{B}}_{#1}}
\newcommand{\RTE}[1]{{\mathrm{RTE}}_{#1}}
\newcommand{\RTF}[1]{{\mathrm{RTF}}_{#1}}
\newcommand{\RTCE}[1]{{\mathrm{RTCE}}_{#1}}
\newcommand{\RTCF}[1]{{\mathrm{RTCF}}_{#1}}
\newcommand{\BDME}[1]{{\mathrm{BDME}}_{#1}}
\newcommand{\BDMF}[1]{{\mathrm{BDMF}}_{#1}}
\newcommand{\BDFM}[1]{{\mathrm{BDFM}}_{#1}}
\newcommand{\NCE}[1]{{\mathrm{NCE}}_{#1}}
\newcommand{\NCF}[1]{{\mathrm{NCF}}_{#1}}
\newcommand{\total}[2]{\frac{\mathrm{d} #1}{\mathrm{d} #2}}
\newcommand{\pp}[2]{\frac{\partial #1}{\partial #2}}
\newcommand{\Hcurl}{H(\operatorname{curl})}
\newcommand{\Hdiv}{H(\operatorname{div})}
\begin{document}

\author[1,2,3,*]{Andrew~T.~T.~McRae}
\author[4]{Gheorghe-Teodor~Bercea}
\author[4,2]{Lawrence~Mitchell}
\author[2,4]{David~A.~Ham}
\author[2]{Colin~J.~Cotter}
\affil[1]{The Grantham Institute, Imperial College London, London, SW7 2AZ, UK}
\affil[2]{Department of Mathematics, Imperial College London, London, SW7 2AZ, UK}
\affil[3]{Department of Mathematical Sciences, University of Bath, Bath, BA2 7AY, UK}
\affil[4]{Department of Computing, Imperial College London, London, SW7 2AZ, UK}
\affil[*]{Correspondence to: \texttt{a.t.t.mcrae@bath.ac.uk}}
\title{Automated generation and symbolic manipulation of tensor product
finite elements}
\date{}
\maketitle

\begin{abstract}
  We describe and implement a symbolic algebra for scalar and vector-valued
  finite elements, enabling the computer generation of elements with tensor
  product structure on quadrilateral, hexahedral and triangular prismatic
  cells. The algebra is implemented as an extension to the domain-specific
  language UFL, the Unified Form Language. This allows users to construct
  many finite element spaces beyond those supported by existing software
  packages. We have made corresponding extensions to FIAT, the FInite
  element Automatic Tabulator, to enable numerical tabulation of such
  spaces. This tabulation is consequently used during the automatic
  generation of low-level code that carries out local assembly operations,
  within the wider context of solving finite element problems posed over
  such function spaces. We have done this work within the code-generation
  pipeline of the software package Firedrake; we make use of the full
  Firedrake package to present numerical examples.
\end{abstract}

\section{Introduction}
\label{sec:intro}

Many different areas of science benefit from the ability to generate
approximate numerical solutions to partial differential equations. In
the past decade, there has been increasing use of software packages and
libraries that automate fundamental operations. The FEniCS
Project~\citep{logg2012automated} is especially notable for allowing the
user to express discretisations of PDEs, based on the finite element
method, in UFL~\citep{alnaes2014unified, alnaes2012ufl} -- a concise,
high-level language embedded in Python. Corresponding efficient
low-level code is automatically generated by FFC, the FEniCS Form
Compiler~\citep{kirby2006compiler, logg2012ffc}, making use of
FIAT~\citep{kirby2004fiat, kirby2012fiat}. These local ``kernels'' are
executed on each cell\footnote{Note on terminology: throughout this
paper, we use the term `cell' to denote the geometric component of the
mesh; we reserve the term `finite element' to denote the space of
functions on a cell and supplementary information related to global
continuity.} in the mesh, and the resulting global systems of equations
can be solved using a number of third-party libraries.

There are multiple advantages to having the discretisation represented
symbolically within a high-level language. The user can write down
complicated expressions concisely without being encumbered by low-level
implementation details. Suitable optimisations can then be applied
automatically during the generation of low-level code; this would be a
tedious process to replicate by hand on each new expression. Such
transformations have previously been implemented in
FFC~\citep{oelgaard2010optimizations, kirby2006compiler}. In this paper,
we extend this high-level approach by introducing a user-facing symbolic representation of tensor product finite elements. Firstly, this enables
the construction of a wide range of finite element spaces, particularly
scalar- and vector-valued identifications of finite element differential
forms. Secondly, while we have not done this at present, the symbolic
representation of a tensor product finite element may be exploited to
automatically generate optimal-complexity algorithms via a
sum-factorisation approach.

Firedrake is an alternative software package to FEniCS which presents a
similar -- in many cases, identical -- interface. Like FEniCS, Firedrake
automatically generates low-level C kernels from high-level UFL
expressions. However, the execution of these kernels over the mesh is
performed in a fundamentally different way; this led to significant
performance increases, relative to FEniCS 1.5, across a range of
problems~\citep{rathgeber2015firedrake}.
As well as the high-level representation of finite element operations
embedded in Python, Firedrake and FEniCS have other attractive features.
They support a wide range of arbitrary-order finite element families,
which are of use to numerical analysts proposing novel discretisations
of PDEs. They also make use of third-party libraries, notably
PETSc~\citep{petsc-user-ref}, exposing a wide range of solvers and
preconditioners for efficient solution of linear systems.

A limitation of Firedrake and FEniCS has been the lack of support for
anything other than fully unstructured meshes with simplicial cells:
intervals, triangles or tetrahedra. There are good reasons why a user
may wish to use a mesh of non-simplicial cells. Our main motivation is
geophysical simulations, which are governed by highly anisotropic
equations in which gravity plays an important role. In addition, they
often require high aspect-ratio domains: the vertical height of the
domain may be several orders of magnitude smaller than the horizontal
width. These domains admit a decomposition which has an unstructured
horizontal `base mesh' but with regular vertical layers -- we will refer
to this as an \emph{extruded} mesh. The cells in such a mesh are not
simplices but instead have a product structure. In two dimensions this
leads to quadrilateral cells; in three dimensions, triangular prisms or
hexahedra. From a mathematical viewpoint, the vertical alignment of
cells minimises difficulties associated with the anisotropy of the
governing equations. From a computational viewpoint, the vertical
structure can be exploited to improve performance compared to a fully
unstructured mesh.

On such cells, we will focus on producing finite elements that can be
expressed as (sums of) products of existing finite elements. This covers
many, though not all, of the common finite element spaces on product
cells. We pay special attention to element families relevant to finite
element exterior calculus, a mathematical framework that leads to stable
mixed finite element discretisations of partial differential
equations~\citep{arnold2006finite, arnold2010finite, arnold2014finite}.
This paper therefore describes some of the extensions to the Firedrake
code-generation pipeline to enable the solution of finite element
problems on cells which are products of simplices. These enable the
automated generation of low-level kernels representing finite element
operations on such cells. We remark that, due to our geophysical
motivations, Firedrake has complete support for extruded meshes whose
unstructured base mesh is built from simplices or quadrilaterals. At the
time of writing, however, it does not support fully unstructured
prismatic or hexahedral meshes.

Many, though not all, of the finite elements we can now construct
already have implementations in other finite element libraries.
deal.II~\citep{dealII82} contains both scalar-valued tensor product
finite elements and the vector-valued Raviart-Thomas and Nédélec
elements of the first kind~\citep{raviart1977mixed, nedelec1980mixed},
which can be constructed using tensor products. However, deal.II only
supports quadrilateral and hexahedral cells and has no support for
simplices or triangular prisms. DUNE PDELab~\citep{bastian2010generic} contains
low-order Raviart-Thomas elements on quadrilaterals and hexahedra, but
only supports scalar-valued elements on triangular prisms.
Nektar++~\citep{cantwell2015nektar} uses tensor-product elements
extensively and supports a wide range of geometric cells, but is
restricted to scalar-valued finite elements. MFEM~\citep{mfem-library}
supports Raviart-Thomas and Nédélec elements of the first kind, though
it has no support for triangular prisms.
NGSolve~\citep{schoberl2014implementation,schoberl2005high} contains many,
possibly all, of the exterior-calculus-inspired tensor-product elements
that we can create on triangular prisms and hexahedra. However, it does
not support elements such as the Nédélec element of the second
kind~\citep{nedelec1986new} on these cells, which do not fit into the
exterior calculus framework.

This paper is structured as follows: in \autoref{sec:prelim}, we provide
the mathematical details of product finite elements. In
\autoref{sec:impl}, we describe the software extensions that allow such
elements to be represented and numerically tabulated. In
\autoref{sec:numex}, we present numerical experiments that make use of
these elements, within Firedrake. Finally, in \autoref{sec:limext} and
\autoref{sec:conc}, we give some limitations of our implementation and
other closing remarks.

\subsection{Summary of contributions}
\begin{itemize}
\item The description and implementation of a symbolic algebra on
existing scalar- and vector-valued finite elements. This allows for the
creation of scalar-valued continuous and discontinuous tensor product
elements, and vector-valued curl- and div-conforming tensor product
elements in two and three dimensions.
\item Certain vector-valued finite elements on quadrilaterals, triangular
prisms and hexahedra are completely unavailable in other major packages,
and some elements we create have no previously published implementation.
\item The tensor product element structure is captured symbolically at
runtime. Although we do not take advantage of this at present, this
could later be exploited to automate the generation of low-complexity
algorithms through sum-factorisation and similar techniques.
\end{itemize}

\section{Mathematical preliminaries}
\label{sec:prelim}

This section is structured as follows: in \autoref{ssec:ciarlet}, we give
the definition of a finite element that we work with. In
\autoref{ssec:sum-elt}, we briefly define the sum of finite elements. In
\autoref{ssec:ss-cont}, we discuss finite element spaces in terms of their
inter-cell continuity. In \autoref{ssec:tp-fem} and
\autoref{ssec:feec-spaces-product}, which form the main part of this
section, we define the product of finite elements and state how these
products can be manipulated and combined to produce elements compatible
with finite element exterior calculus. Up to this point, our exposition
uses the language of scalar and vector fields as our existing software
infrastructure uses scalars and vectors and we believe this makes the
paper accessible to a wider audience. However, we end this section with
\autoref{ssec:tp-complex-forms}, which briefly re-states
\autoref{ssec:tp-fem} and \autoref{ssec:feec-spaces-product} in terms of
differential forms. These provide a far more natural setting for the
underlying operations.

\subsection{Definition of a finite element}
\label{ssec:ciarlet}

We will follow \citet{ciarlet1978finite} in defining a
\emph{finite element} to be a triple ($K$,~$P$,~$N$) where
\begin{itemize}
  \item $K$ is a bounded domain in~$\mathbb{R}^n$, to be interpreted as
  a generic \emph{reference cell} on which all calculations are
  performed,
  \item $P$ is a finite-dimensional space of continuous functions
  on~$K$, typically some subspace of polynomials,
  \item $N = \{n_1, \ldots, n_{\dim P}\}$ is a basis for the dual
  space~$P'$ -- the space of linear functionals on~$P$ -- where the
  elements of the set~$N$ are called \emph{nodes}.
\end{itemize}

Let~$\Omega$ be a compact domain which is decomposed into a finite
number of non-overlapping cells. Assume that we wish to find an
approximate solution to some partial differential equation, posed
in~$\Omega$, using the finite element method. A \emph{finite element}
together with a given decomposition of~$\Omega$ produce a \emph{finite
element space}.

A finite element space is a finite-dimensional function space
on~$\Omega$. There are essentially two things that need to be specified
to characterise a finite element space: the manner in which a function
may vary within a single cell, and the amount of continuity a function
must have between neighbouring cells.

The former is related to~$P$; more details are given in
\autoref{ssec:piola}. A basis for~$P$ is therefore very useful in
implementations of the finite element method. Often, this is a
\emph{nodal basis} in which each of the basis
functions~$\Phi_1, \ldots, \Phi_{\dim P}$ vanish when acted on by all
but one node:
\begin{equation}
n_i(\Phi_j) = \delta_{ij}.
\end{equation}

Basis functions from different cells can be combined into basis
functions for the finite element space on~$\Omega$. The inter-cell
continuity of these basis functions is related to the choice of
nodes,~$N$. This is the core topic of \autoref{ssec:ss-cont}.

\subsection{Sum of finite elements}
\label{ssec:sum-elt}

Suppose we have finite elements~$U = (K, P_A, N_A)$
and~$V = (K, P_B, N_B)$, which are defined over the same reference
cell~$K$. If the intersection of~$P_A$ and~$P_B$ is trivial, we can
define the \emph{direct sum}~${U \oplus V}$ to be the
finite element~$(K, P, N)$, where
\begin{align}
P &\vcentcolon = P_A \oplus P_B \equiv \{f_A + f_B \mid f_A \in P_A, f_B \in P_B\} \\
N &\vcentcolon = N_A \cup N_B.
\end{align}

\subsection{Sobolev spaces, inter-cell continuity, and Piola transforms}
\label{ssec:ss-cont}

Finite element spaces are a finite-dimensional subspace of some larger
Sobolev space, depending on the degree of continuity of functions
between neighbouring cells. We will consider finite element spaces in
$H^1$, $\Hcurl$, $\Hdiv$ and $L^2$.

A brief remark: it is clear that these Sobolev spaces have some trivial
inclusion relations -- $H^1$ is a subspace of $L^2$, $\Hdiv$ and $\Hcurl$
are both subspaces of $[L^2]^d$, where $d$ is the spatial dimension, and
$[H^1]^d$ is a subspace of both $\Hdiv$ and $\Hcurl$. However, in what
follows, when we make casual statements such
as~$V \subset \Hdiv$, it is \emph{implied} that~$V \not\subset [H^1]^d$,
i.e., we have made the strongest statement possible. In particular, we
will use~$L^2$ to denote a total absence of continuity between cells.

\subsubsection{Geometric decomposition of nodes}

The set of nodes~$N$, from the definition in \autoref{ssec:ciarlet}, are
used to enforce the continuity requirements on the `global' finite
element space. This is done by associating nodes with \emph{topological}
entities of~$K$ -- vertices, facets, and so on. When multiple cells
in~$\Omega$ share a topological entity, the cells must agree on the
value of any degree of freedom associated with that entity. This leads
to coupling between any cells that share the entity. The association of
nodes with topological entities is crucial in determining the continuity
of finite element spaces -- this is sometimes called the
\emph{geometric decomposition} of nodes.

For~$H^1$ elements, functions are fully continuous between cells, and
must therefore be single-valued on vertices, edges and facets. Nodes are
firstly associated with \emph{vertices}. If necessary, additional nodes
are associated with \emph{edges}, then with \emph{facets}, then with the
\emph{interior} of the reference cell.

For~$\Hcurl$ elements, which are intrinsically vector-valued, functions
must have continuous tangential component between cells. The
component(s) of the function tangential to edges and facets must
therefore be single-valued. Nodes are firstly associated with
\emph{edges} until the tangential component is specified uniquely. If
necessary, additional nodes are associated with \emph{facets}, then with
the \emph{interior} of the reference cell.

For~$\Hdiv$ elements, which are also intrinsically vector-valued,
functions must have continuous normal component between cells. The
component of the function normal to facets must therefore be
single-valued. Nodes are firstly associated with \emph{facets}. If
necessary, additional nodes are associated with the \emph{interior} of
the cell.

$L^2$ elements have no continuity requirements. Typically, all nodes are
associated with the \emph{interior} of the cell; this does not lead to
any continuity constraints.

\subsubsection{Piola transforms}
\label{ssec:piola}

For functions to have the desired continuity on the global mesh, they
may need to undergo an appropriate \emph{mapping} from reference to
physical space. Let~$\vec{X}$ represent coordinates on the reference
cell, and~$\vec{x}$ represent coordinates on the physical cell; for each
physical cell there is some map~$\vec{x} = g(\vec{X})$.

For $H^1$ or~$L^2$ functions, no explicit mapping is needed.
Let~$\hat{f}(\vec{X})$ be a function defined over the reference cell.
The corresponding function~$f(\vec{x})$ defined over the physical cell
is then
\begin{equation}
\label{eq:identity-map}
f(\vec{x}) = \hat{f} \circ g^{-1}(\vec{x}).
\end{equation}
We will refer to this as the \emph{identity} mapping.

However, if we wish to have continuity of the normal or tangential
component of the vector field in physical space;
Eq.~\eqref{eq:identity-map} does not suffice. $\Hdiv$~and~$\Hcurl$
elements therefore use \emph{Piola transforms} to map functions from
reference space to physical space. We will use~$J$ to
denote~$Dg(\vec{X})$, the Jacobian of the coordinate transformation.
$\Hdiv$ functions are mapped using the contravariant Piola transform,
which preserves normal components:
\begin{equation}
\label{eq:cont-Piola}
\vec{f}(\vec{x}) = \frac{1}{\det J} J \hat{\vec{f}} \circ g^{-1}(\vec{x}),
\end{equation}
while~$\Hcurl$ functions are mapped using the covariant Piola transform,
which preserves tangential components:
\begin{equation}
\label{eq:co-Piola}
\vec{f}(\vec{x}) = J^{-T}\hat{\vec{f}} \circ g^{-1}(\vec{x}).
\end{equation}

\subsection{Product finite elements}
\label{ssec:tp-fem}

In this section, we discuss how to take the product of a pair of finite
elements and how this product element may be manipulated to give
different types of inter-cell continuity. We will label our constituent
elements~$U$ and~$V$, where $U \vcentcolon = (K_A, P_A, N_A)$ and
$V \vcentcolon = (K_B, P_B, N_B)$ following the notation
of \autoref{ssec:ciarlet}. We begin with the definition of the product
reference cell, which is straightforward. However, the spaces of
functions and the associated nodes are intimately related, hence the
discussion of these is interleaved.

\subsubsection{Product cells}
\label{sssec:prod-cell}

Given reference cells $K_A \subset \mathbb{R}^n$ and
$K_B \subset \mathbb{R}^m$, the reference product cell $K_A \times K_B$
can be defined straightforwardly as follows:
\begin{equation}
\label{eq:prod-cell}
K_A \times K_B \vcentcolon = \left\{ (x_1, \ldots, x_{n+m}) \in \mathbb{R}^{n+m} \mid (x_1, \ldots, x_n) \in K_A, (x_{n+1}, \ldots, x_{n+m}) \in K_B \right\}.
\end{equation}

The topological entities of $K_A \times K_B$ correspond to products of
topological entities of~$K_A$ and~$K_B$. If we label the entities of a
reference cell (in $\mathbb{R}^n$, say) by their dimension, so that $0$
corresponds to vertices, $1$ to edges, \ldots, $n-1$ to facets and $n$
to the cell, the entities of $K_A \times K_B$ can be labelled as follows:
\begin{description}
  \item[(0, 0):] vertices of $K_A \times K_B$ -- the product of a vertex of $K_A$ with a vertex of $K_B$
  \item[(1, 0):] edges of $K_A \times K_B$ -- the product of an edge of $K_A$ with a vertex of $K_B$
  \item[(0, 1):] edges of $K_A \times K_B$ -- the product of a vertex of $K_A$ with an edge of $K_B$\\ 
  \vdots
  \item[(n-1, m):] facets of $K_A \times K_B$ -- the product of a facet of $K_A$ with the cell of $K_B$
  \item[(n, m-1):] facets of $K_A \times K_B$ -- the product of the cell of $K_A$ with a facet of $K_B$
  \item[(n, m):] cell of $K_A \times K_B$ -- the product of the cell of $K_A$ with the cell of $K_B$
\end{description}
It is important to distinguish between different types of entities, even
those with the same dimension. For example, if~$K_A$ is a triangle
and~$K_B$ an interval, the~$(2, 0)$ facets of the prism~$K_A \times K_B$
are triangles while the~$(1, 1)$ facets are quadrilaterals.

\subsubsection{Product spaces of functions -- simple elements}
\label{sssec:prod-func-i}

Given spaces of functions~$P_A$ and~$P_B$, the product
space~$P_A \otimes P_B$ can be defined as the span of products of
functions in~$P_A$ and~$P_B$:
\begin{equation}
\label{eq:prod-elt}
P_A \otimes P_B \vcentcolon = \operatorname{span}\left\{ f \cdot g \mid f \in P_A, g \in P_B \right\},
\end{equation}
where the product function $f \cdot g$ is defined so that
\begin{equation}
\label{eq:prod-fn}
(f \cdot g)(x_1, \ldots, x_{n+m}) = f(x_1, \ldots, x_n) \cdot g(x_{n+1}, \ldots, x_{n+m}).
\end{equation}
In the cases we consider explicitly, at least one of~$f$~or~$g$ will be
scalar-valued, so the product on the right-hand side
of Eq.~\eqref{eq:prod-fn} is unambiguous. A basis for $P_A \otimes P_B$ can
be constructed from bases for~$P_A$ and~$P_B$. If $P_A$ and $P_B$ have
nodal bases
\begin{equation}
\label{eq:vw-basis}
\left\{\Phi_1^{(A)}, \Phi_2^{(A)}, \ldots \Phi_N^{(A)} \right\},
\left\{\Phi_1^{(B)}, \Phi_2^{(B)}, \ldots \Phi_M^{(B)} \right\}
\end{equation}
respectively, a nodal basis for $P_A \otimes P_B$ is given by
\begin{equation}
\label{eq:vtimesw-basis}
\left\{\Phi_{i, j},\quad i = 1, \ldots, N, j = 1, \ldots, M \right\},
\end{equation}
where
\begin{equation}
\label{eq:prod-basisfn}
\Phi_{i, j} \vcentcolon = \Phi_i^{(A)}\cdot \Phi_j^{(B)},\quad i = 1, \ldots, N, j = 1, \ldots, M;
\end{equation}
the right-hand side uses the same product as Eq.~\eqref{eq:prod-fn}.

While this already gives plenty of flexibility, there are cases in which
a different, more natural, space can be built by further manipulation
of~$P_A \otimes P_B$. We will return to this after a brief description
of product nodes.

\subsubsection{Product nodes -- geometric decomposition}
\label{sssec:prod-node-i}

Recall that the nodes are a basis for the dual
space~$(P_A \otimes P_B)'$, and that the inter-cell continuity of the
finite element space is related to the association of nodes with
topological entities of the reference cell.

Assuming that we know bases for~$P_A'$ and~$P_B'$, there is a natural
basis for~$(P_A \otimes P_B)'$ which is essentially an outer (tensor) product of
the bases for~$P_A'$ and~$P_B'$. Let~$n_{i, j}$ denote a ``product''
of~$n^{(A)}_i$, the~$i$'th node in~$N_A$, with~$n^{(B)}_j$, the~$j$'th
node in~$N_B$ -- typically the evaluation of some component of the
function. If $n^{(A)}_i$ is associated with an entity of~$K_A$ of
dimension $p$ and $n^{(B)}_j$ is associated with an entity of~$K_B$ of
dimension $q$ then $n_{i, j}$ is associated with an entity
of~$K_A \times K_B$ with label $(p, q)$.

This geometric decomposition of nodes in the product element is used to
motivate further manipulation of~$P_A \otimes P_B$ to produce
a more natural space of functions, particularly in the case of
vector-valued elements.

\subsubsection{Product spaces of functions -- scalar- and vector-valued elements in 2D and 3D}
\label{sssec:prod-func-ii}

In two dimensions, we take the reference cells~$K_A$ and~$K_B$ to be
intervals, so the product cell~$K_A \times K_B$ is two-dimensional.
Finite elements on intervals are scalar-valued and are either in~$H^1$
or~$L^2$. We will consider the creation of two-dimensional elements in
$H^1$, $\Hcurl$, $\Hdiv$ and $L^2$. A summary of the following is given
in \autoref{tbl:2D}.

\vspace{2mm}
\noindent \begin{minipage}{0.74\columnwidth}
$H^1$: The element must have nodes associated with vertices of the
reference product cell. The vertices of the reference product cell are
formed by taking the product of vertices on the intervals. The
constituent elements must therefore have nodes associated with vertices,
so must both be in~$H^1$.
\end{minipage}
\begin{minipage}{0.25\columnwidth}
\centering
\includegraphics[width=0.4\columnwidth]{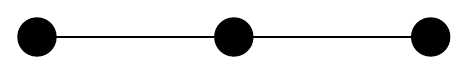} \quad
\includegraphics[width=0.4\columnwidth]{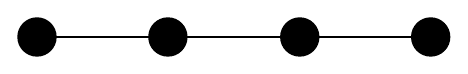} \\
\vspace{3mm}
\includegraphics[width=0.63\columnwidth]{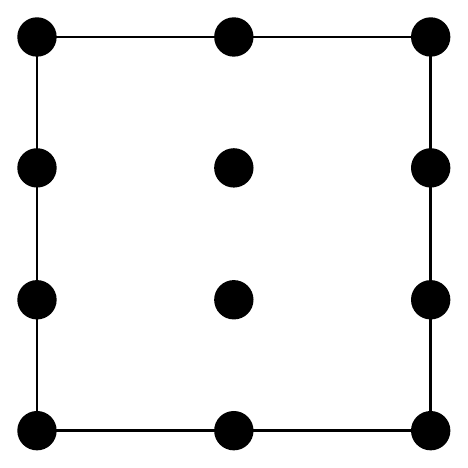}
\end{minipage}

\vspace{2mm}
\noindent \begin{minipage}{0.74\columnwidth}
$\Hcurl$: The element must have nodes associated with edges of the
reference product cell. The edges of the reference product cell are
formed by taking the product of an interval's vertex with an interval's
interior. One of the constituent elements must therefore have nodes
associated with vertices, while the other must only have nodes
associated with the interior. Taking the product of an $H^1$ element
with an $L^2$ element gives a scalar-valued element with nodes on the
$(0, 1)$ facets, for example.

\vspace{2mm}
To create an $\Hcurl$ element, we now multiply this scalar-valued
element by the vector $(0, 1)$ to create a vector-valued finite element
(if we had taken the product of an $L^2$ element with an $H^1$ element,
we would multiply by $(1, 0)$). This gives an element whose tangential
component is continuous across \emph{all} edges (trivially so on two of
the edges). In addition, we must use an appropriate Piola transform when
mapping from reference space into physical space.
\end{minipage}
\begin{minipage}{0.25\columnwidth}
\centering
\includegraphics[width=0.4\columnwidth]{diagrams/h1_int.pdf} \quad
\includegraphics[width=0.4\columnwidth]{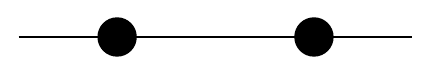} \\
\vspace{3mm}
\includegraphics[width=0.63\columnwidth]{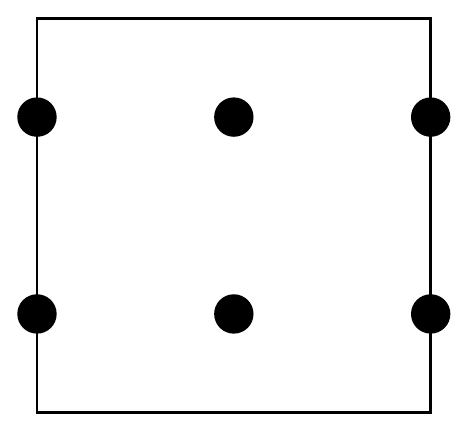} \\
\vspace{3mm}
\includegraphics[width=0.64\columnwidth]{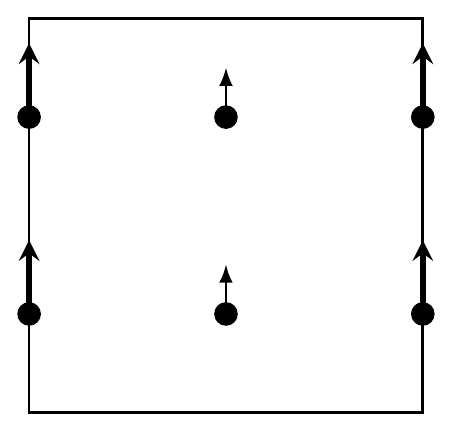}
\end{minipage}

\vspace{2mm}
\noindent \begin{minipage}{0.74\columnwidth}
$\Hdiv$: We create a scalar-valued element in the same way as in the
$\Hcurl$ case, but multiplied by the `other' basis vector (for
$H^1 \times L^2$, we choose $(-1, 0)$ -- the minus sign is
useful for orientation consistency in unstructured quadrilateral meshes;
for $L^2 \times H^1$, $(0, 1)$). This gives
an element whose normal component is continuous across \emph{all} edges,
and again, we must use an appropriate Piola transform when mapping from
reference space into physical space.
\end{minipage}
\begin{minipage}{0.25\columnwidth}
\centering
\includegraphics[width=0.72\columnwidth]{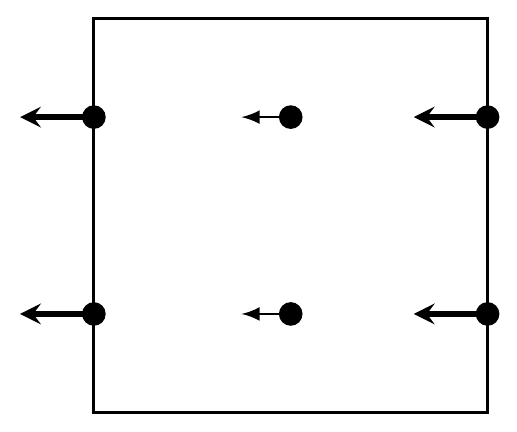}
\end{minipage}

\vspace{2mm}
\noindent
Note that the scalar-valued product elements we produce above are
perfectly legitimate finite elements, and it is not compulsory to form
vector-valued elements from them. Indeed, we use such a scalar-valued
element for the example in \autoref{ssec:ex-grav-wave}. However, the
vector-valued elements are generally more useful and fit naturally
within Finite Element Exterior Calculus, as we will see in
\autoref{ssec:feec-spaces-product}.

\vspace{1mm}
\noindent \begin{minipage}{0.74\columnwidth}
$L^2$: The element must only have nodes associated with interior of the
reference product cell. The constituent elements must therefore only
have nodes associated with their interiors, so must both be in~$L^2$.
\end{minipage}
\begin{minipage}{0.25\columnwidth}
\centering
\includegraphics[width=0.4\columnwidth]{diagrams/l2_int.pdf} \quad
\includegraphics[width=0.4\columnwidth]{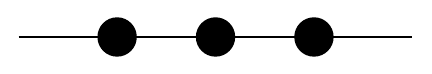} \\
\vspace{3mm}
\includegraphics[width=0.63\columnwidth]{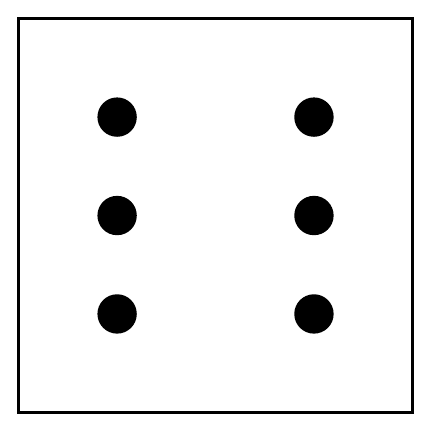}
\end{minipage}
\clearpage

\begin{table}
\footnotesize
\centering
\caption{Summary of 2D product elements\label{tbl:2D}}
\begin{tabu}
{|c|c|c|c|c|}
\hline
Product (1D $\times$ 1D) & Components & Modifier & Result & Mapping\\
\hline
$H^1 \times H^1$ & $f \times g$ & (none) & $fg$ & identity \\
$H^1 \times L^2$ & $f \times g$ & (none) & $fg$ & identity \\
$H^1 \times L^2$ & $f \times g$ & $\Hcurl$ & $(0, fg)$ & covariant Piola \\
$H^1 \times L^2$ & $f \times g$ & $\Hdiv$ & $(-fg, 0)$ & contravariant Piola \\
$L^2 \times H^1$ & $f \times g$ & (none) & $fg$ & identity \\
$L^2 \times H^1$ & $f \times g$ & $\Hcurl$ & $(fg, 0)$ & covariant Piola \\
$L^2 \times H^1$ & $f \times g$ & $\Hdiv$ & $(0, fg)$ & contravariant Piola \\
$L^2 \times L^2$ & $f \times g$ & (none) & $fg$ & identity \\
\hline
\end{tabu}
\end{table}
\vspace{2mm}

In three dimensions, we take~$K_A \subset \mathbb{R}^2$ and~$K_B$ to be
an interval, so the product cell~$K_A \times K_B$ is three-dimensional.
Finite elements on a 2D reference cell may be in~$H^1$, $\Hcurl$,
$\Hdiv$ or~$L^2$. Elements on a 1D reference cell may be in~$H^1$
or~$L^2$. We will consider the creation of three-dimensional elements in
$H^1$, $\Hcurl$, $\Hdiv$ and $L^2$. A summary of the following is given
in \autoref{tbl:3D}.
\vspace{2mm}

\emph{Note: In the following pictures, we have taken the two-dimensional
cell to be a triangle. However, the discussion is equally valid for
quadrilaterals.}

\noindent \begin{minipage}{0.74\columnwidth}
$H^1$: As in the two-dimensional case, this is formed by taking the
product of two $H^1$ elements.
\end{minipage}
\begin{minipage}{0.25\columnwidth}
\centering
\includegraphics[width=0.25\columnwidth]{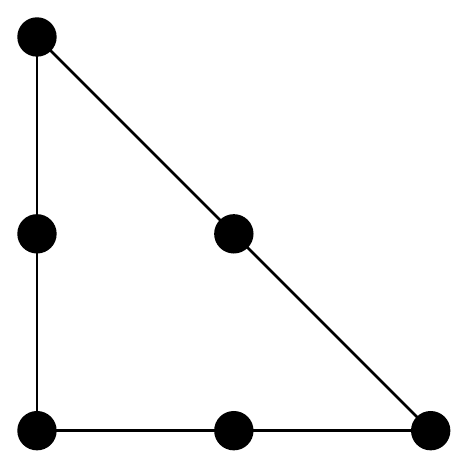} \quad
\includegraphics[width=0.4\columnwidth]{diagrams/h1_int.pdf} \\
\vspace{1mm}
\includegraphics[width=0.5\columnwidth]{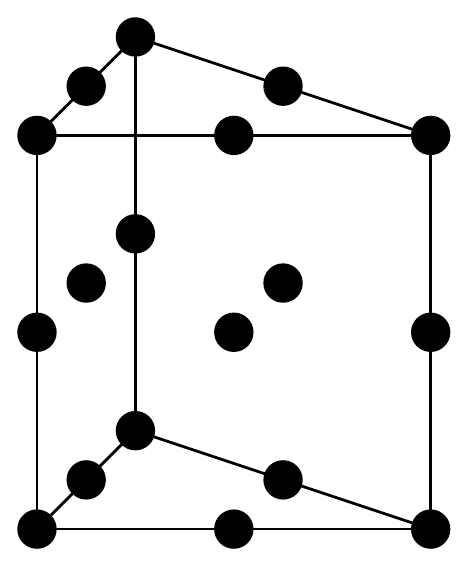}
\end{minipage}

\noindent \begin{minipage}{0.74\columnwidth}
$\Hcurl$: The element must again have nodes associated with edges of the
reference product cell. There are two distinct ways of forming such an
element, and in both cases a suitable Piola transform must be used to
map functions from reference to physical space.

\vspace{2mm}
Taking the product of an $H^1$ two-dimensional element with an
$L^2$ one-dimensional element produces a scalar-valued element with
nodes on $(0, 1)$ edges. If we multiply this by the vector $(0, 0, 1)$,
this results in an element whose tangential component is continuous on
all edges and faces.
\end{minipage}
\begin{minipage}{0.25\columnwidth}
\centering
\includegraphics[width=0.25\columnwidth]{diagrams/h1_tri.pdf} \quad
\includegraphics[width=0.4\columnwidth]{diagrams/l2_int.pdf} \\
\vspace{1mm}
\includegraphics[width=0.5\columnwidth]{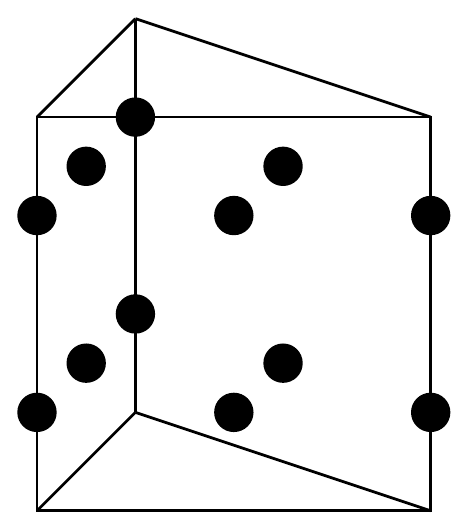} \\
\vspace{1mm}
\includegraphics[width=0.5\columnwidth]{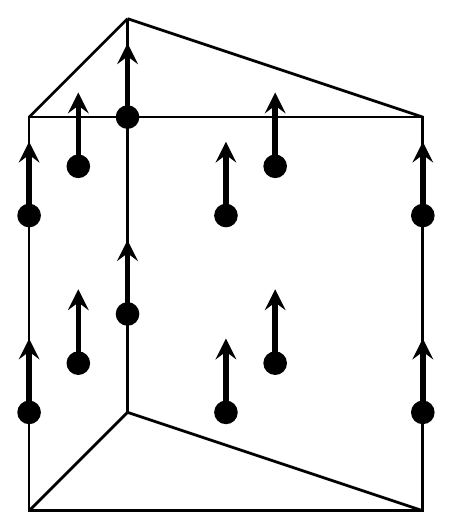}
\end{minipage}

\noindent \begin{minipage}{0.74\columnwidth}
Alternatively, one may take the product of an $\Hdiv$ or $\Hcurl$
two-dimensional element with an $H^1$ one-dimensional element. This
produces a vector-valued element with nodes on $(1, 0)$ edges. The
product naturally takes values in~$\mathbb{R}^2$, since the
two-dimensional element is vector-valued and the one-dimensional element
is scalar-valued. However, an $\Hcurl$ element in three dimensions must
take values in~$\mathbb{R}^3$. If the two-dimensional element is in
$\Hcurl$, it is enough to interpret the product as the first two
components of a three-dimensional vector. If the two-dimensional element
is in $\Hdiv$, the two-dimensional product must be \emph{rotated} by 90
degrees before being transformed into a three-dimensional vector.
\end{minipage}
\begin{minipage}{0.25\columnwidth}
\centering
\includegraphics[width=0.3\columnwidth]{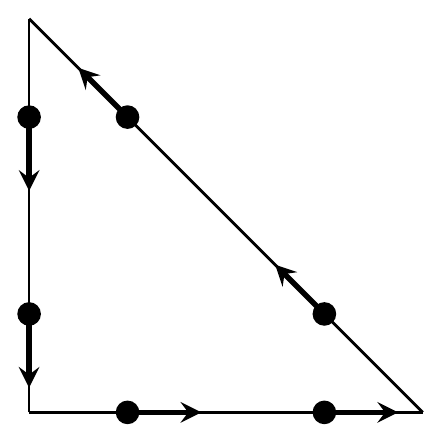} \quad
\includegraphics[width=0.4\columnwidth]{diagrams/h1_int.pdf} \\
\vspace{1mm}
\includegraphics[width=0.5\columnwidth]{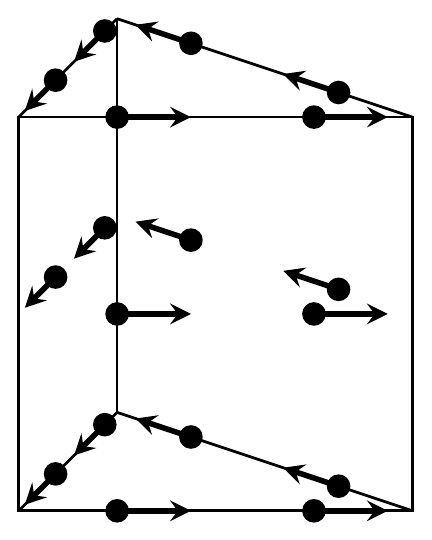}
\end{minipage}

\noindent \begin{minipage}{0.74\columnwidth}
$\Hdiv$: The element must have nodes associated with facets of the
reference product cell. As with $\Hcurl$, there are two distinct ways of
forming such an element, and suitable Piola transforms must again be
used.

\vspace{2mm}
Taking the product of an $L^2$ two-dimensional element with an $H^1$
one-dimensional element gives a scalar-valued element with nodes on
$(2, 0)$ facets. Multiplying this by $(0, 0, 1)$ produces an element
whose normal component is continuous across all facets.
\end{minipage}
\begin{minipage}{0.25\columnwidth}
\centering
\includegraphics[width=0.25\columnwidth]{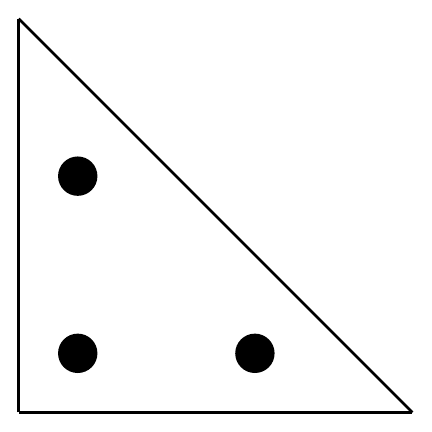} \quad
\includegraphics[width=0.4\columnwidth]{diagrams/h1_int.pdf} \\
\vspace{1mm}
\includegraphics[width=0.5\columnwidth]{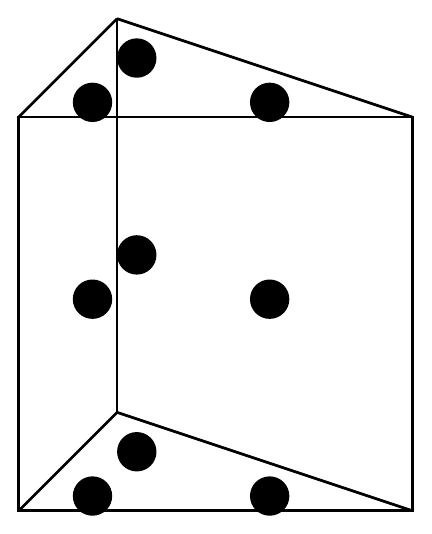}
\vspace{1mm}
\includegraphics[width=0.5\columnwidth]{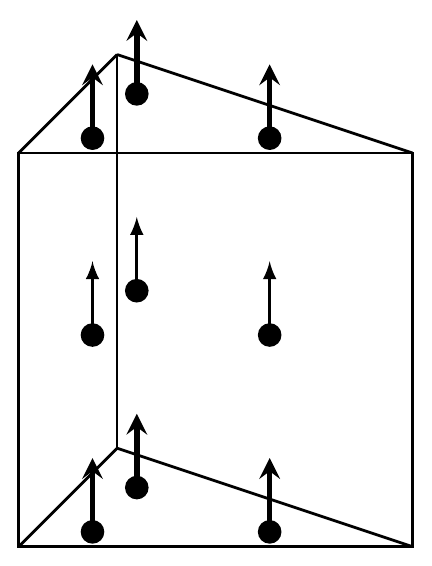}
\end{minipage}

\noindent \begin{minipage}{0.74\columnwidth}
Taking the product of an $\Hdiv$ or $\Hcurl$ two-dimensional element
with an $L^2$ one-dimensional element gives a vector-valued element with
nodes on $(1, 1)$ facets. Again, the product naturally takes values
in~$\mathbb{R}^2$. If the two-dimensional element is in $\Hdiv$, it is
enough to interpret the product as the first two components of a
three-dimensional vector-valued element whose third component vanishes.
If the two-dimensional element is in $\Hcurl$, the product must be
rotated by 90 degrees before transforming.
\end{minipage}
\begin{minipage}{0.25\columnwidth}
\centering
\includegraphics[width=0.35\columnwidth]{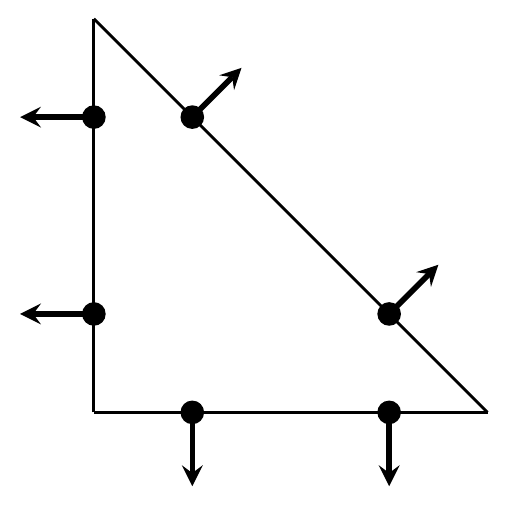} \quad
\includegraphics[width=0.3\columnwidth]{diagrams/l2_int.pdf} \\
\vspace{1mm}
\includegraphics[width=0.55\columnwidth]{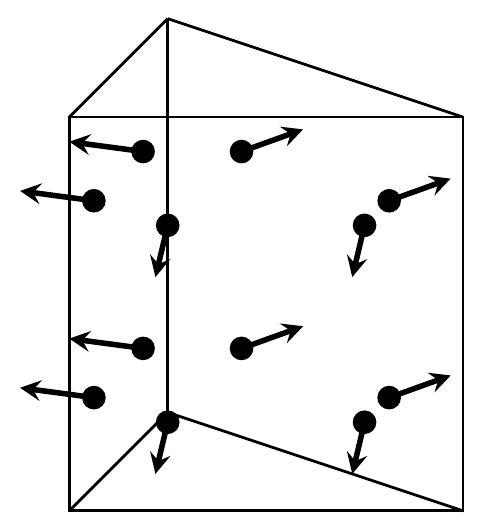}
\end{minipage}

\noindent \begin{minipage}{0.74\columnwidth}
$L^2$: As in the 2D case, both constituent elements must be
in $L^2$.
\end{minipage}
\begin{minipage}{0.25\columnwidth}
\centering
\includegraphics[width=0.25\columnwidth]{diagrams/l2_tri.pdf} \quad
\includegraphics[width=0.4\columnwidth]{diagrams/l2_int.pdf} \\
\vspace{1mm}
\includegraphics[width=0.5\columnwidth]{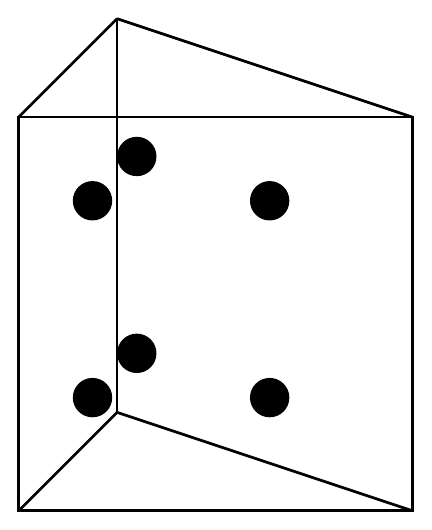}
\end{minipage}

\begin{table}
\footnotesize
\centering
\caption{Summary of 3D product elements\label{tbl:3D}}
\begin{tabu}{|c|c|c|c|c|}
\hline
Product (2D $\times$ 1D) & Components & Modifier & Result & Mapping\\
\hline
$H^1 \times H^1$ & $f \times g$ & (none) & $fg$ & identity \\
$H^1 \times L^2$ & $f \times g$ & (none) & $fg$ & identity \\
$H^1 \times L^2$ & $f \times g$ & $\Hcurl$ & $(0, 0, fg)$ & covariant Piola \\
$\Hcurl \times H^1$ & $(f_x, f_y) \times g$ & (none) & $(f_x g, f_y g)^\dagger$ & * \\
$\Hcurl \times H^1$ & $(f_x, f_y) \times g$ & $\Hcurl$ & $(f_x g, f_y g, 0)$ & covariant Piola \\
$\Hdiv \times H^1$ & $(f_x, f_y) \times g$ & (none) & $(f_x g, f_y g)^\dagger$ & * \\
$\Hdiv \times H^1$ & $(f_x, f_y) \times g$ & $\Hcurl$ & $(-f_y g, f_x g, 0)$ & covariant Piola \\
$\Hcurl \times L^2$ & $(f_x, f_y) \times g$ & (none) & $(f_x g, f_y g)^\dagger$ & * \\
$\Hcurl \times L^2$ & $(f_x, f_y) \times g$ & $\Hdiv$ & $(f_y g, -f_x g, 0)$ & contravariant Piola \\
$\Hdiv \times L^2$ & $(f_x, f_y) \times g$ & (none) & $(f_x g, f_y g)^\dagger$ & * \\
$\Hdiv \times L^2$ & $(f_x, f_y) \times g$ & $\Hdiv$ & $(f_x g, f_y g, 0)$ & contravariant Piola \\
$L^2 \times H^1$ & $f \times g$ & (none) & $fg$ & identity \\
$L^2 \times H^1$ & $f \times g$ & $\Hdiv$ & $(0, 0, fg)$ & contravariant Piola \\
$L^2 \times L^2$ & $f \times g$ & (none) & $fg$ & identity \\
\hline
\end{tabu}

The elements marked with $^\dagger$ are of little practical use; they
are 2-vector valued but are defined over three-dimensional domains. No
mapping has been given for these elements; the Piola transformations
from a 3D cell require all three components to be defined.
\end{table}

\subsubsection{Consequences for implementation}
\label{ssec:consequences}
The previous subsections motivate the implementation of several
mathematical operations on finite elements. We will need an operator
that takes the product of two existing elements; we call this
\texttt{TensorProductElement}. This will generate a new element whose
reference cell is the product of the reference cells of the constituent
elements, as described in \autoref{sssec:prod-cell}. It will also
construct the product space of functions~$P_A \otimes P_B$, as described
in \autoref{sssec:prod-func-i}, but with no extra manipulation (e.g.\ 
expanding into a vector-valued space). The basis for $P_A \otimes P_B$
is as defined in Eqs.~\eqref{eq:vtimesw-basis}
and~\eqref{eq:prod-basisfn}. The nodes are topologically associated with
topological entities of the reference cell as described in
\autoref{sssec:prod-node-i}.

To construct the more complicated vector-valued finite elements, we
introduce additional operators \texttt{HCurl} and \texttt{HDiv} which
form a vector-valued $\Hcurl$ or $\Hdiv$ element from an existing
\texttt{TensorProductElement}. This will modify the product space as
described in \autoref{sssec:prod-func-ii} by manipulating the existing
product into a vector of the correct dimension (after rotation, if
applicable), and setting an appropriate Piola transform. We will also
need an operator that creates the sum of finite elements; this already
exists in UFL under the name \texttt{EnrichedElement}, and is
represented by \texttt{+}.

\subsection{Product finite elements within finite element exterior calculus}
\label{ssec:feec-spaces-product}

The work of \citet{arnold2006finite, arnold2010finite} on
finite element exterior calculus provides principles for obtaining
stable mixed finite element discretisations on a domain consisting of
simplicial cells: intervals, triangles, tetrahedra, and
higher-dimensional analogues. In full generality, this involves de~Rham
complexes of polynomial-valued finite element differential forms linked
by the exterior derivative operator. In 1, 2 and 3 dimensions,
differential forms can be naturally identified with scalar and vector
fields, while the exterior derivative can be interpreted as a standard
differential operator such as grad, curl, or div. The vector-valued
element spaces only have partial continuity between cells: they are
in~$\Hcurl$ or~$\Hdiv$, which have been discussed already. The element
spaces themselves were, however, already well-known in the existing
finite element literature for their use in solving mixed formulations of
the Poisson equation and problems of a similar nature.

\citet{arnold2014finite} generalises finite element
exterior calculus to cells which can be expressed as geometric products
of simplices. It also describes a specific complex of finite element
spaces on hexahedra (and, implicitly, quadrilaterals). When these
differential forms are identified with scalar- and vector-valued
functions, they correspond to the scalar-valued~$\Q{r}$ and
its discontinuous counterpart $\DQ{r}$, and various well-known
vector-valued spaces as introduced in \citet{brezzi1985two},
\citet{nedelec1980mixed} and \citet{nedelec1986new}. Within finite element
exterior calculus, there are element spaces which cannot be expressed as
a tensor product of spaces on simplices -- see, for example,
\citet{arnold2014finite2} -- but we are not considering such
spaces in this paper.

Finite element exterior calculus makes use of de~Rham complexes of
finite element spaces. In one dimension, the complex takes the form
\begin{equation}
\label{eq:recap-complex-1d}
U_0 \stackrel{\total{}{x}}{\longrightarrow} U_1,
\end{equation}
where $U_0 \subset H^1$ and $U_1 \subset L^2$. In two dimensions, there
are two types of complex, arising due to two possible identifications of
differential 1-forms with vector fields:
\begin{equation}
\label{eq:recap-complex-2d-div}
U_0 \stackrel{\nabla^\perp}{\longrightarrow} U_1 \stackrel{\nabla\cdot}{\longrightarrow} U_2,
\end{equation}
where $U_0 \subset H^1$, $U_1 \subset \Hdiv$, and $U_2 \subset L^2$, and
\begin{equation}
\label{eq:recap-complex-2d-curl}
U_0 \stackrel{\nabla}{\longrightarrow} U_1 \stackrel{\nabla^\perp\cdot}{\longrightarrow} U_2,
\end{equation}
where $U_0 \subset H^1$, $U_1 \subset \Hcurl$, and $U_2 \subset L^2$. In
three dimensions, the complex takes the form
\begin{equation}
\label{eq:recap-complex-3d}
U_0 \stackrel{\nabla}{\longrightarrow} U_1 \stackrel{\nabla\times}{\longrightarrow}U_2 \stackrel{\nabla\cdot}{\longrightarrow} U_3,
\end{equation}
where $U_0 \subset H^1$, $U_1 \subset \Hcurl$, $U_2 \subset \Hdiv$, and
$U_3 \subset L^2$.

Given an existing two-dimensional complex $(U_0, U_1, U_2)$ and a
one-dimensional complex $(V_0, V_1)$, we can generate a product complex
on the three-dimensional product cell:
\begin{equation}
\label{eq:prod-complex-3d}
W_0 \stackrel{\nabla}{\longrightarrow} W_1 \stackrel{\nabla\times}{\longrightarrow} W_2 \stackrel{\nabla\cdot}{\longrightarrow} W_3,
\end{equation}
where
\begin{align}
W_0 &\vcentcolon = U_0 \otimes V_0, \\
W_1 &\vcentcolon = \texttt{HCurl}(U_0 \otimes V_1) \oplus \texttt{HCurl}(U_1 \otimes V_0), \\
W_2 &\vcentcolon = \texttt{HDiv}(U_1 \otimes V_1) \oplus \texttt{HDiv}(U_2 \otimes V_0), \\
W_3 &\vcentcolon = U_2 \otimes V_1,
\end{align}
with
$W_0 \subset H^1, W_1 \subset \Hcurl, W_2 \subset \Hdiv, W_3 \subset L^2$
(compare the complex given in Eq.~\eqref{eq:recap-complex-3d}). The
vector-valued spaces are direct sums of `product' spaces that have
been modified by the \texttt{HCurl} or \texttt{HDiv} operator.

Similarly, taking the product of two one-dimensional complexes produces
a product complex on the two-dimensional product cell in which the
vector-valued space is in \emph{either} $\Hdiv$ or $\Hcurl$.

\subsection{Product complexes using differential forms}
\label{ssec:tp-complex-forms}

This section summarises \citet{arnold2014finite} by
restating the results of \autoref{ssec:tp-fem} and
\autoref{ssec:feec-spaces-product} in the language of differential forms,
which can be considered a generalisation of scalar and vector fields.

In three dimensions, 0-forms and 3-forms are identified with scalar
fields, while 1-forms and 2-forms are identified with vector fields. In
two dimensions, 0-forms and 2-forms are identified with scalar fields.
1-forms are identified with vector fields, but this can be done in two
different ways since 1-forms and ($n$-1)-forms coincide. This results in
two possible vector fields, which differ by a 90-degree rotation. In one
dimension, both 0-forms and 1-forms are conventionally identified with
scalar fields.

Let~$K_A\subset\mathbb{R}^n$, $K_B\subset\mathbb{R}^m$ be domains.
Suppose we are given de~Rham subcomplexes on~$K_A$ and~$K_B$,
\begin{align}
\label{eq:cpx}
U_0 \stackrel{\mathrm{d}}{\longrightarrow} U_1 \stackrel{\mathrm{d}}{\longrightarrow} \cdots \stackrel{\mathrm{d}}{\longrightarrow} U_n, \qquad
V_0 \stackrel{\mathrm{d}}{\longrightarrow} V_1 \stackrel{\mathrm{d}}{\longrightarrow} \cdots \stackrel{\mathrm{d}}{\longrightarrow} V_m,
\end{align}
where each $U_k$ is a space of (polynomial) differential $k$-forms on
$K_A$ and each $V_k$ is a space of differential $k$-forms on $K_B$.
The product of these complexes is a de~Rham subcomplex on
$K_A \times K_B$:
\begin{equation}
\label{eq:cpx-prdef}
(U \otimes V)_0 \stackrel{\mathrm{d}}{\longrightarrow} (U \otimes V)_1 \stackrel{\mathrm{d}}{\longrightarrow} \cdots \stackrel{\mathrm{d}}{\longrightarrow} (U \otimes V)_{n+m},
\end{equation}
where, for $k = 0, 1, \ldots, n+m$,
\begin{equation}
\label{eq:cpx-otdef}
(U \otimes V)_k \vcentcolon = \bigoplus_{i+j=k} (U_i \otimes V_j).
\end{equation}

Note that~$(U \otimes V)_k$ is a space of (polynomial) $k$-forms
on~$K_A \otimes K_B$, and can hence be interpreted as a scalar or vector
field in 2 or 3 spatial dimensions. It can be easily verified
that the definitions in Eqs.~\eqref{eq:cpx-prdef} and~\eqref{eq:cpx-otdef}
gives rise to Eq.~\eqref{eq:prod-complex-3d} in three dimensions, for
example. The discussion
in \autoref{sssec:prod-func-i} and \autoref{sssec:prod-func-ii} on the
product of function spaces can be summarised by the definition
of~$\otimes$ on the right-hand side of Eq.~\eqref{eq:cpx-otdef}, along with
the definition of the standard wedge product of differential forms.
It is clear that much of the apparent complexity of the \texttt{HDiv}
and \texttt{HCurl} operators introduced in \autoref{ssec:tp-fem} arises
from working with scalars and vectors rather than introducing
differential forms!

\section{Implementation}
\label{sec:impl}
The symbolic operations on finite elements, derived in the previous
section, have been implemented within
Firedrake~\citep{rathgeber2015firedrake, rathgeber2014productive}.
Firedrake is an ``\emph{automated system for the portable solution of
partial differential equations using the finite element method}".
Firedrake has several dependencies. Some of these are components of the
FEniCS Project~\citep{logg2012automated}:
\begin{description}
\item[FIAT] FInite element Automatic
Tabulator~\citep{kirby2004fiat, kirby2012fiat}, for the construction
and tabulation of finite element basis functions
\item[UFL] Unified Form
Language~\citep{alnaes2014unified, alnaes2012ufl}, a domain-specific
language for the specification of finite element variational forms
\end{description}
Firedrake also relies on PyOP2~\citep{rathgeber2012pyop2} and
COFFEE~\citep{luporini2015crossloop}.

The changes required to effect the generation of product elements were
largely confined to FIAT and UFL, while support for integration over
product cells is included in Firedrake's form compiler. We begin this
section with
more detailed expositions on FIAT and UFL. We discuss the implementation
of product finite elements in \autoref{ssec:impl-prod-elt}. We talk about
the resulting algebraic structure in \autoref{ssec:algebraic-structure}.
We finish by discussing the new integration regions, in
\autoref{ssec:impl-int-reg}.

\subsection{FIAT}
\label{ssec:FIAT}

This component is responsible for computing finite element basis
functions for a wide range of finite element families. To do this, it
works with an abstraction based on Ciarlet's definition of a finite
element, as given in \autoref{ssec:ciarlet}. The reference cell~$K$ is
defined using a set of vertices, with higher-dimensional geometrical
objects defined as sets of vertices. The polynomial space~$P$ is defined
implicitly through a \emph{prime} basis: typically an orthonormal set
of polynomials, such as (on triangles) a Dubiner basis, which can be
stably evaluated to high polynomial order. The set of nodes~$N$ is
also defined; this implies the existence of a \emph{nodal} basis
for~$P$, as explained previously.

The nodal basis, which is important in calculations, can be expressed
as linear combinations of prime basis functions. This is done
automatically by FIAT; details are given in \citet{kirby2004fiat}.
The main method of interacting with FIAT is by
requesting the tabulated values of the nodal basis functions at a set of
points inside~$K$ -- typically a set of quadrature points. FIAT also
stores the geometric decomposition of nodes relative to the topological
entities of~$K$.

\subsection{UFL}
\label{ssec:UFL}

This component is a domain-specific language, embedded in Python, for
representing weak formulations of partial differential equations. It is
centred around expressing multilinear forms: maps from the product of
some set of function spaces~$\{V_j\}^\rho_{j=1}$ into the real numbers
which are linear in each argument, where $\rho$ is 0, 1 or 2.
Additionally, the form may be parameterised over one or more
\emph{coefficient functions}, and is not necessarily linear in these.
The form may include derivatives of functions, and the language has
extensive support for matrix algebra operations.

We can assume that the function spaces are finite element spaces; in
UFL, these are represented by the \texttt{FiniteElement} class. This
requires three pieces of information: the element family, the geometric
cell, and the polynomial degree. A limited amount of symbolic
manipulation on \texttt{FiniteElement} objects could already be done:
the UFL \texttt{EnrichedElement} class is used to represent the~$\oplus$
operator discussed in \autoref{ssec:sum-elt}.

\subsection{Implementation of product finite elements}
\label{ssec:impl-prod-elt}

To implement product finite elements, additions to UFL and FIAT were
required. The UFL changes are purely symbolic and allow the new elements
to be represented. The FIAT changes allow the new elements (and
derivatives thereof) to be numerically tabulated at specified points in
the reference cell.

As discussed in \autoref{ssec:consequences}, we implemented several new
element classes in UFL. The existing UFL \texttt{FiniteElement} classes
has two essential properties: the \texttt{degree} and the
\texttt{value\_shape}. The \texttt{degree} is the maximal degree of any
polynomial basis function -- this allows determination of an appropriate
quadrature rule. The \texttt{value\_shape} represents whether the
element is scalar-valued or vector-valued and, if applicable, the
dimension of the vector in \emph{physical} space. This allows suitable
code to be generated when doing vector and tensor operations.

For \texttt{TensorProductElement}s, we define the \texttt{degree} to be a
tuple; the basis functions are products of polynomials in distinct sets
of variables. It is therefore advantageous to store the polynomial
degrees separately for later use with a product quadrature rule. The
\texttt{value\_shape} is defined according to the definition in
\autoref{sssec:prod-func-i} for the product of functions. For
\texttt{HCurl} and \texttt{HDiv} elements, the \texttt{degree} is
identical to the \texttt{degree} of the underlying
\texttt{TensorProductElement}. The \texttt{value\_shape} needs to be
modified: in physical space, these vector-valued elements have dimension
equal to the dimension of the physical space.

The secondary role of FIAT is to store a representation of the geometric
decomposition of nodes. For product elements, the generation of this was
described in \autoref{sssec:prod-node-i}. The primary role is to tabulate
finite element basis functions, and derivatives thereof, at specified
points in the reference cell. The \texttt{tabulate} method of a FIAT
finite element takes two arguments: the maximal \texttt{order} of
derivatives to tabulate, and the set of \texttt{points}.

Let $\Phi_{i, j}(x, y, z) \vcentcolon= \Phi_i^{(A)}(x, y)\Phi_j^{(B)}(z)$
be some product element basis function; we will assume that this is
scalar-valued to ease the exposition. Suppose we need to tabulate the
$x$-derivative of this at some specified point $(x_0, y_0, z_0)$.
Clearly
\begin{equation}
\pp{\Phi_{i, j}}{x}(x_0, y_0, z_0) = \pp{\Phi_i^{(A)}}{x}(x_0, y_0)\Phi_j^{(B)}(z_0).
\end{equation}
In other words, the value can be obtained from tabulating (derivatives
of) basis functions of the constituent elements at appropriate points.
It is clear that this extends to other combinations of derivatives, as
well as to components of vector-valued basis functions. Further
modifications to the tabulation for curl- or div-conforming vector
elements are relatively simple, as detailed in
\autoref{sssec:prod-func-ii}.

\subsection{Algebraic structure}
\label{ssec:algebraic-structure}

The extensions described in \autoref{ssec:impl-prod-elt} enable
sophisticated manipulation of finite elements within UFL. For example,
consider the following complex on triangles, highlighted
by \citet{cotter2012mixed} as being relevant for numerical weather
prediction:
\begin{equation}
\P{2} \oplus \B{3} \stackrel{\nabla^\perp}{\longrightarrow} \BDFM{2} \stackrel{\nabla\cdot}{\longrightarrow} \DP{1}.
\end{equation}
Here, $\P{2} \oplus \B{3}$ denotes the space of quadratic polynomials
enriched by a cubic `bubble' function, $\BDFM{2}$ represents a member of
the vector-valued Brezzi--Douglas--Fortin--Marini element
family~\citep{brezzi1991mixed} in~$\Hdiv$, and $\DP{1}$ represents the
space of discontinuous, piecewise-linear functions. Suppose we wish to
take the product of this with some complex on intervals, such as
\begin{equation}
\P{2} \stackrel{\total{}{x}}{\longrightarrow} \DP{1}.
\end{equation}
This generates a complex on triangular prisms:
\begin{equation}
W_0 \stackrel{\nabla}{\longrightarrow} W_1 \stackrel{\nabla\times}{\longrightarrow} W_2 \stackrel{\nabla\cdot}{\longrightarrow} W_3,
\end{equation}
where
\begin{align}
W_0 &\vcentcolon = (\P{2}^\triangle \oplus \B{3}^\triangle) \otimes \P{2}, \\
W_1 &\vcentcolon = \texttt{HCurl}((\P{2}^\triangle \oplus \B{3}^\triangle) \otimes \DP{1}) \oplus \texttt{HCurl}(\BDFM{2}^\triangle \otimes \P{2}), \\
W_2 &\vcentcolon = \texttt{HDiv}(\BDFM{2}^\triangle \otimes \DP{1}) \oplus \texttt{HDiv}(\DP{1}^\triangle \otimes \P{2}), \\
W_3 &\vcentcolon = \DP{1}^\triangle \otimes \DP{1};
\end{align}
we have marked the elements on triangles by $^\triangle$ for clarity.
Following our extensions to UFL, the product complex may be constructed
as shown in Listing~\ref{lst:prod-comp-ex}. Some of these elements are
used in the example in \autoref{ssec:ex-grav-wave}.

\begin{lstlisting}[float, frame=single, language={[firedrake]{python}}, label=lst:prod-comp-ex, caption=Construction of a complicated product complex in UFL]
U0_0 = FiniteElement("P", triangle, 2)
U0_1 = FiniteElement("B", triangle, 3)
U0 = EnrichedElement(U0_0, U0_1)
U1 = FiniteElement("BDFM", triangle, 2)
U2 = FiniteElement("DP", triangle, 1)

V0 = FiniteElement("P", interval, 1)
V1 = FiniteElement("DP", interval, 0)

W0 = TensorProductElement(U0, V0)
W1_h = TensorProductElement(U1, V0)
W1_v = TensorProductElement(U0, V1)
W1 = EnrichedElement(HCurl(W1_h), HCurl(W1_v))
W2_h = TensorProductElement(U1, V1)
W2_v = TensorProductElement(U2, V0)
W2 = EnrichedElement(HDiv(W2_h), HDiv(W2_v))
W3 = TensorProductElement(U2, V1)
\end{lstlisting}

\subsection{Support for new integration regions}
\label{ssec:impl-int-reg}

On simplicial meshes, Firedrake supports three types of integrals:
integrals over cells, integrals over exterior facets and integrals over
interior facets. Integrals over exterior facets are typically used to
apply boundary conditions weakly, while integrals over interior facets
are used to couple neighbouring cells when discontinuous function spaces
are present. The implementation of the different types of integral is
quite elegant: the only difference between integrating a function over
the interior of the cell and over a single facet is the choice of
quadrature points and quadrature weights.  Note that Firedrake assumes
that the mesh is conforming -- hanging nodes are not currently supported.

On product cells, all entities can be considered as a product of
entities on the constituent cells. We can therefore construct product
quadrature rules, making use of existing quadrature rules for
constituent cells and facets thereof. In addition, we split the facet
integrals into separate integrals over `vertical' and `horizontal' facets.
This is natural when executing a computational kernel over an extruded
unstructured mesh, and may be useful in geophysical contexts where
horizontal and vertical motions may be treated differently.

\section{Numerical examples}
\label{sec:numex}

In this section, we give several examples to demonstrate the correctness
of our implementation. Quantitative analysis is performed where
possible, e.g. demonstration of convergence to a known solution at
expected order with increasing mesh resolution. Tests are performed in
both two and three spatial dimensions. We make use of Firedrake's
\texttt{ExtrudedMesh} functionality. In two dimensions, the cells are
quadrilaterals, usually squares. In three dimensions, we use triangular
prisms, though we can also build elements on hexahedra.

When referring to standard finite element spaces, we follow the
convention in which the number refers to the degree of the minimal
complete polynomial space containing the element, not the maximal
complete polynomial space contained by the element. Thus, an element
containing some, but not all, linear polynomials is numbered~$1$,
rather than~$0$. This is the convention used by UFL, and is also
justified from the perspective of finite element exterior calculus.

\subsection{Vector Laplacian (3D)}
\label{ssec:vlap}
We seek a solution to
\begin{equation}
\label{eq:v-lap-strong}
-\nabla(\nabla\cdot\vec{u}) + \nabla\times(\nabla\times\vec{u}) = \vec{f}
\end{equation}
in a domain $\Omega$, with boundary conditions
\begin{align}
\label{eq:v-lap-bcs-1}
\vec{u} \cdot \vec{n} &= 0, \\
\label{eq:v-lap-bcs-2}
(\nabla \times \vec{u}) \times \vec{n} &= 0
\end{align}
on $\partial{\Omega}$, where $\vec{n}$ is the outward normal. A naïve
discretisation can lead to spurious solutions, especially on non-convex
domains, but an accurate discretisation can be obtained by introducing
an auxiliary variable (see, for example, \citet{arnold2010finite}):
\begin{align}
\label{eq:v-lap-strong-alt}
\sigma &= -\nabla\cdot\vec{u}, \\
\nabla \sigma + \nabla\times(\nabla\times\vec{u}) &= \vec{f}.
\end{align}

Let $V_0 \subset H^1$, $V_1 \subset \Hcurl$ be finite element spaces. A
suitable weak formulation is: find $\sigma \in V_0$, $\vec{u} \in V_1$
such that
\begin{align}
\label{eq:v-lap-weak}
\langle \tau, \sigma \rangle - \langle \nabla\tau, \vec{u} \rangle &= 0, \\
\langle \vec{v}, \nabla \sigma \rangle + \langle \nabla\times\vec{v}, \nabla\times\vec{u} \rangle &= \langle \vec{v}, \vec{f} \rangle,
\end{align}
for all $\tau \in V_0, \vec{v} \in V_1$, where we have used angled
brackets to denote the standard $L^2$ inner product. The boundary
conditions have been implicitly applied, in a weak sense, through
neglecting the surface terms when integrating by parts.

We take $\Omega$ to be the unit cube $[0, 1]^3$. Let $k$, $l$ and $m$ be
arbitrary. Then
\begin{equation}
\label{eq:v-lap-f-3d}
\vec{f} = \pi^2\begin{pmatrix}
(k^2 + l^2)\sin(k\pi x)\cos(l\pi y)\\
(l^2 + m^2)\sin(l\pi y)\cos(m\pi z)\\
(k^2 + m^2)\sin(m\pi z)\cos(k\pi x)
\end{pmatrix}
\end{equation}
produces the solution
\begin{equation}
\label{eq:v-lap-u-3d}
\vec{u} = \begin{pmatrix}
\sin(k\pi x)\cos(l\pi y)\\
\sin(l\pi y)\cos(m\pi z)\\
\sin(m\pi z)\cos(k\pi x)
\end{pmatrix},
\end{equation}
which satisfies the boundary conditions.

\begin{figure}
\centering
\includegraphics[width=0.75\columnwidth]{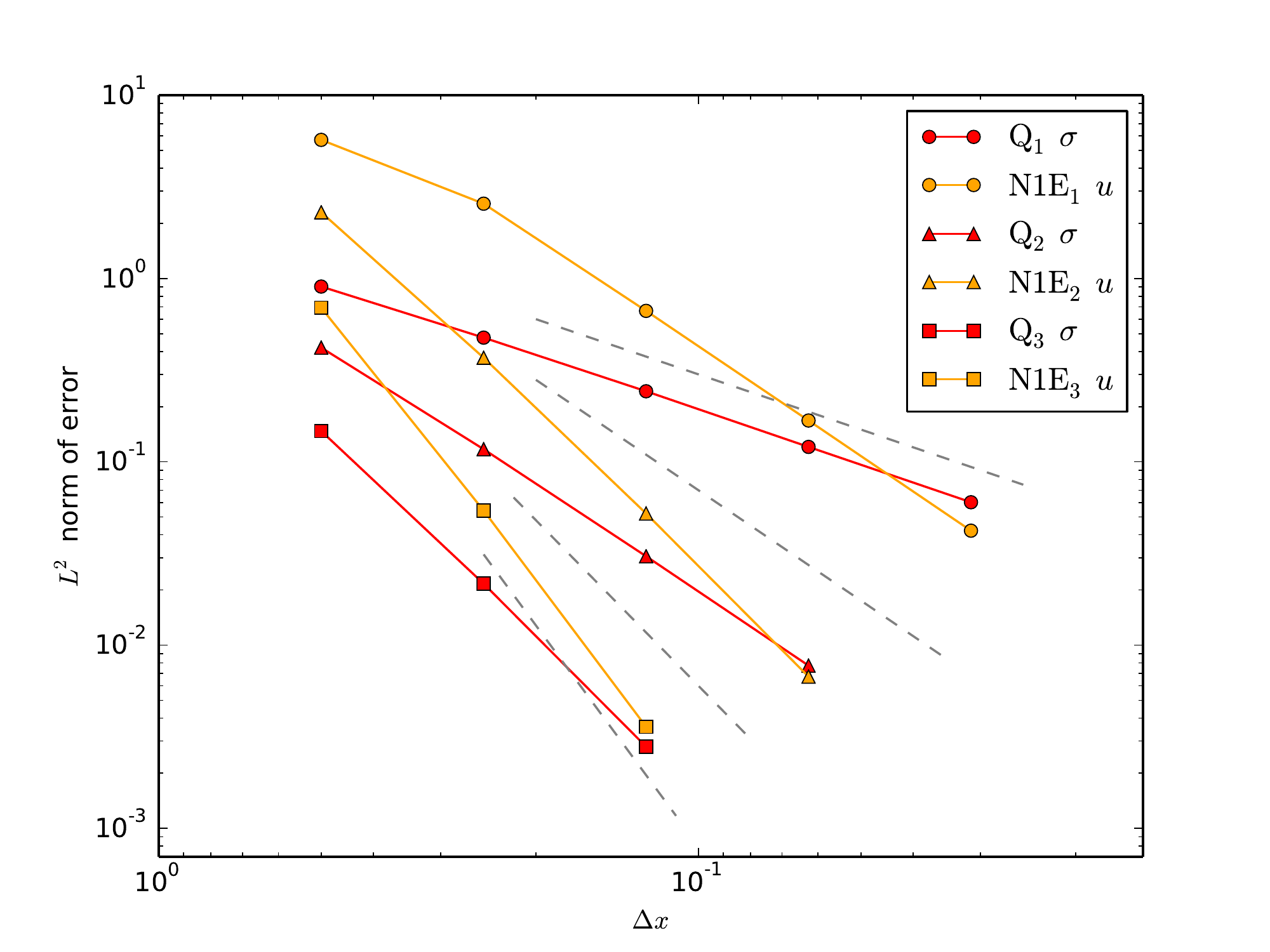}
\caption{The $L^2$ error between the computed and `analytic' solution
is plotted against~$\Delta x$ for the 3D problem described
in \autoref{ssec:vlap}. The dotted lines are
proportional to~$\Delta x^n$, for~$n$ from 1 to 4, and are merely to
aid comprehension.}
\label{fig:ex1_plot3d}
\end{figure}

To discretise this problem, we subdivide $\Omega$ into triangular
prisms whose base is a right-angled triangle with short sides of
length~$\Delta x$ and whose height is~$\Delta x$. We use
the~$\Q{r}$ prism element for the~$H^1$ space, and the degree-$r$
Nédélec prism element of the first kind for the~$\Hcurl$ space, for~$r$
from 1 to 3. We take $k$, $l$ and $m$ to be 1, 2 and 3, respectively.
We approximate~$\vec{f}$ by interpolating the analytic
expression onto a vector-valued function in~$\Q{r+1}$.
The~$L^2$ errors between the calculated and `analytic' solutions for
varying $\Delta x$ are plotted in \autoref{fig:ex1_plot3d}. This is
done for both~$\vec{u}$ and~$\sigma$; the so-called analytic solutions
are approximations which are formed by interpolating the genuine
analytic solution onto nodes of~$\Q{r+1}$.

\subsection{Gravity wave (3D)}
\label{ssec:ex-grav-wave}

A simple model of atmospheric flow is given by
\begin{equation}
\label{eq:grav-wave}
\pp{\vec{u}}{t} = -\nabla p + b \hat{z},\quad
\pp{b}{t} = -N^2 \vec{u} \cdot \hat{z},\quad
\pp{p}{t} = -c^2 \nabla \cdot \vec{u},
\end{equation}
along with the boundary condition $\vec{u}\cdot\vec{n} = 0$,
where~$\vec{n}$ is a unit normal vector. The prognostic variables are
the velocity,~$\vec{u}$, the pressure perturbation,~$p$, and the
buoyancy perturbation,~$b$. The scalars $N$ and~$c$ are (dimensional)
constants, while~$\hat{z}$ represents a unit vector opposite to the
direction of gravity. These equations are a reduction of, for example,
(17)--(21) from \citet{skamarock1994efficiency}, in which we have
neglected the constant background velocity and the Coriolis term, and
rescaled~$\theta$ by~$\theta_0 / g$ to produce~$b$.

Given some three-dimensional product complex as in
Eq.~\eqref{eq:prod-complex-3d}, we seek a solution with
$\vec{u} \in W_2^0$, $b \in W_2^v$ and $p \in W_3$. $W_2^0$ is the
subspace of~$W_2$ whose normal component vanishes on the boundary of the
domain. $W_2^v$ denotes the ``vertical" part of~$W_2$: if we write~$W_2$
as a sum of two product elements~$\texttt{HDiv}(U_1 \otimes V_1)$
and~$\texttt{HDiv}(U_2 \otimes V_0)$ then $W_2^v$ is the
\emph{scalar-valued} product~$U_2 \otimes V_0$, as was constructed in
Listing~\ref{lst:prod-comp-ex}. This combination of finite element
spaces for~$\vec{u}$ and~$b$ is analogous to the Charney--Phillips
staggering of variables in the vertical
direction \citep{charney1953numerical}.

A semi-discrete form of~\eqref{eq:grav-wave} is the following:
find~$\vec{u} \in W_2^0$, $b \in W_2^v$, $p \in W_3$ such that for all
$\vec{w} \in W_2^0$, $\gamma \in W^v_2$,
$\phi \in W_3$
\begin{align}
\label{eq:grav-wave-disc-begin}
\left\langle \vec{w}, \pp{\vec{u}}{t}\right\rangle - \left\langle \nabla\cdot\vec{w}, p\right\rangle - \left\langle \vec{w}, b \hat{z}\right\rangle &= 0 \\
\left\langle \gamma, \pp{b}{t}\right\rangle + N^2\left\langle \gamma, \vec{u} \cdot \hat{z}\right\rangle &= 0 \\
\label{eq:grav-wave-disc-end}
\left\langle \phi, \pp{p}{t}\right\rangle + c^2 \left\langle \phi, \nabla \cdot \vec{u}\right\rangle &= 0.
\end{align}

It can be easily verified that the original
equations,~\eqref{eq:grav-wave}, together with the boundary condition
lead to conservation of the energy perturbation
\begin{equation}
\label{eq:energy}
\ints{\Omega} \frac{1}{2}|\vec{u}|^2 + \frac{1}{2N^2}b^2 + \frac{1}{2c^2}p^2 \dx.
\end{equation}
The three terms can be interpreted as kinetic energy (KE), potential
energy (PE) and internal energy (IE), respectively. The
semi-discretisation given in
Eqs.~\eqref{eq:grav-wave-disc-begin}--\eqref{eq:grav-wave-disc-end}
also conserves this energy. If we discretise in time using the implicit
midpoint rule, which preserves quadratic
invariants \citep{leimkuhler2005simulating} then the fully discrete
system will conserve energy as well.

We take the domain to be a spherical shell centred at the origin. Its
inner radius, $a$, is approximately 6371km, and its thickness, $H$, is
10km. The domain is divided into triangular prism cells with
side-lengths of the order of 1000km and height 1km. We take
$N = 10^{-2}$s$^{-1}$ and $c = 300$ms$^{-1}$. The simulation starts at
rest with a buoyancy perturbation and a vertically balancing pressure
field given by
\begin{equation}
\label{eq:grav-ics-3d}
b = \frac{\sin(\pi(|\vec{x}| - a)/H)}{1+z^2/L^2}, \qquad
p = -\frac{H}{\pi} \frac{\cos(\pi(|\vec{x}| - a)/H)}{1+z^2/L^2};
\end{equation}
$L$ is a horizontal length-scale, which we take to be 500km. We use a
timestep of 1920s, and run for a total of 480,000s.

\begin{figure}
\centering
\includegraphics[width=0.75\columnwidth]{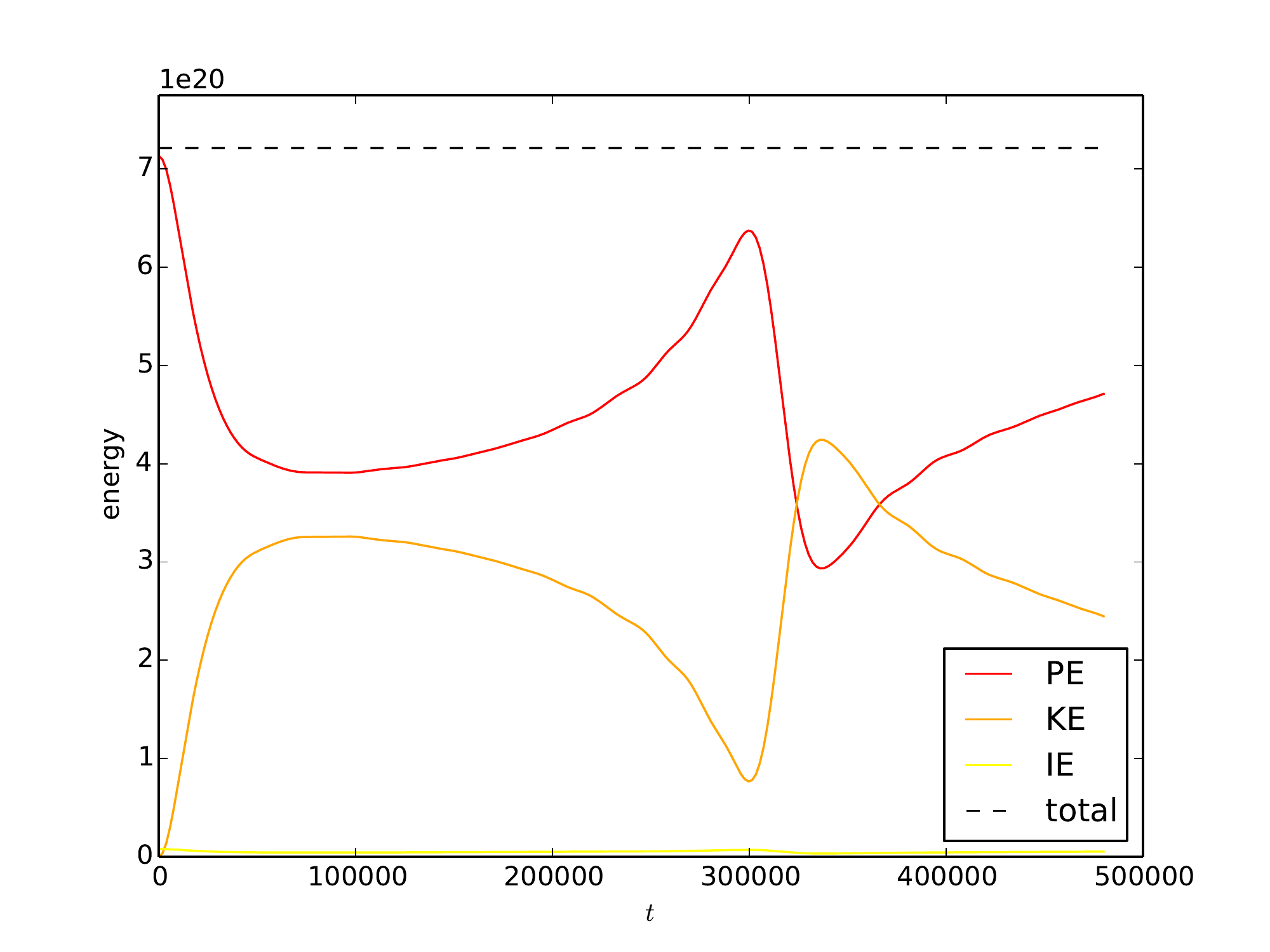}
\caption{Evolution of energy for the simulation described
in \autoref{ssec:ex-grav-wave}. The components are the potential energy, PE,
the kinetic energy, KE, and the internal energy, IE. The choice of
spatial and temporal discretisations leads to exact conservation of
total energy up to solver tolerances; this is indeed observed. The event
at approximately $t$ = 320,000s corresponds to the zonally-symmetric
gravity wave reaching the poles of the spherical domain.}
\label{fig:ex2_plot3d}
\end{figure}

To discretise this problem, we use the product elements formed from the
$\BDFM{2}$ complex on triangles and the $\P{2}$--$\DP{1}$ complex on
intervals; these were constructed in \autoref{ssec:algebraic-structure}.
The initial conditions are interpolated into the buoyancy and pressure
fields. The energy is calculated at every time step; the results are
plotted in \autoref{fig:ex2_plot3d}. The total energy is conserved to
roughly one part in $1.4 \times 10^8$, which is comparable to the linear
solver tolerances.

\subsection{DG advection (2D)}
\label{sec:dgadv-2d}

The advection of a scalar field~$q$ by a known divergence-free velocity
field~$\vec{u_0}$ can be described by the equation
\begin{equation}
\label{eq:dg-cont}
\frac{\partial{q}}{\partial{t}} + \nabla\cdot(\vec{u_0}q) = 0.
\end{equation}
If~$q$ is in a discontinuous function space,~$V$, a suitable weak
formulation is
\begin{equation}
\label{eq:dg-disc-final}
\left\langle \phi, \frac{\partial{q}}{\partial{t}} \right\rangle = \left\langle \nabla \phi, q \vec{u_0} \right\rangle\, -\, \int_{\Gamma_\mathrm{ext}} \! \phi \widetilde{q} \vec{u_0} \cdot \vec{n} \, \mathrm{d} s\, -\, \int_{\Gamma_\mathrm{int}} \! \llbracket\phi\rrbracket \widetilde{q} \vec{u_0} \cdot \vec{n} \, \mathrm{d} S,
\end{equation}
for all $\phi \in V$, where the integrals on the right hand side are
over \emph{exterior} and \emph{interior} mesh facets, with $\mathrm{d}s$
and $\mathrm{d}S$ appropriate integration measures. $\vec{n}$ is the
appropriately-oriented normal vector, $\widetilde{q}$ represents the
\emph{upwind} value of $q$, while $\llbracket\phi\rrbracket$
represents the jump in $\phi$. We assume that, on parts of the boundary
corresponding to inflow, $\widetilde{q} = 0$. This example will therefore
demonstrate the ability to integrate over interior and exterior mesh facets.

We discretise Eq.~\eqref{eq:dg-disc-final} in time using the third-order
three-stage strong-stability-preserving Runge-Kutta scheme given
in \citet{shu1988efficient}. We take $\Omega$ to be the unit square
$[0, 1]^2$. Our initial condition will be a cosine hill
\begin{equation}
\label{eq:dgadv-2d-q-ic}
q = \begin{cases} \frac{1}{2}\left(1 + \cos\left(\pi\frac{|\vec{x} - \vec{x_0}|}{r_0}\right)\right), &|\vec{x} - \vec{x_0}| < r_0 \\
0, &\mathrm{otherwise},
\end{cases}
\end{equation}
with radius~$r_0 = 0.15$, centred at~$\vec{x_0} = (0.25, 0.5)$. The
prescribed velocity field is
\begin{equation}
\label{eq:dgadv-2d-u0}
\vec{u_0}(\vec{x}, t) = \cos\left(\frac{\pi t}{T}\right)\begin{pmatrix}
\sin(\pi x)^2\sin(2\pi y)\\
-\sin(\pi y)^2\sin(2\pi x)
\end{pmatrix},
\end{equation}
as in \citet{leveque1996high}. This gives a reversing, swirling flow
field which vanishes on the boundaries of~$\Omega$. The initial
condition should be recovered at~$t = T$.
Following \citep{leveque1996high}, we take~$T=\frac{3}{2}$.

\begin{figure}
\centering
\includegraphics[width=0.75\columnwidth]{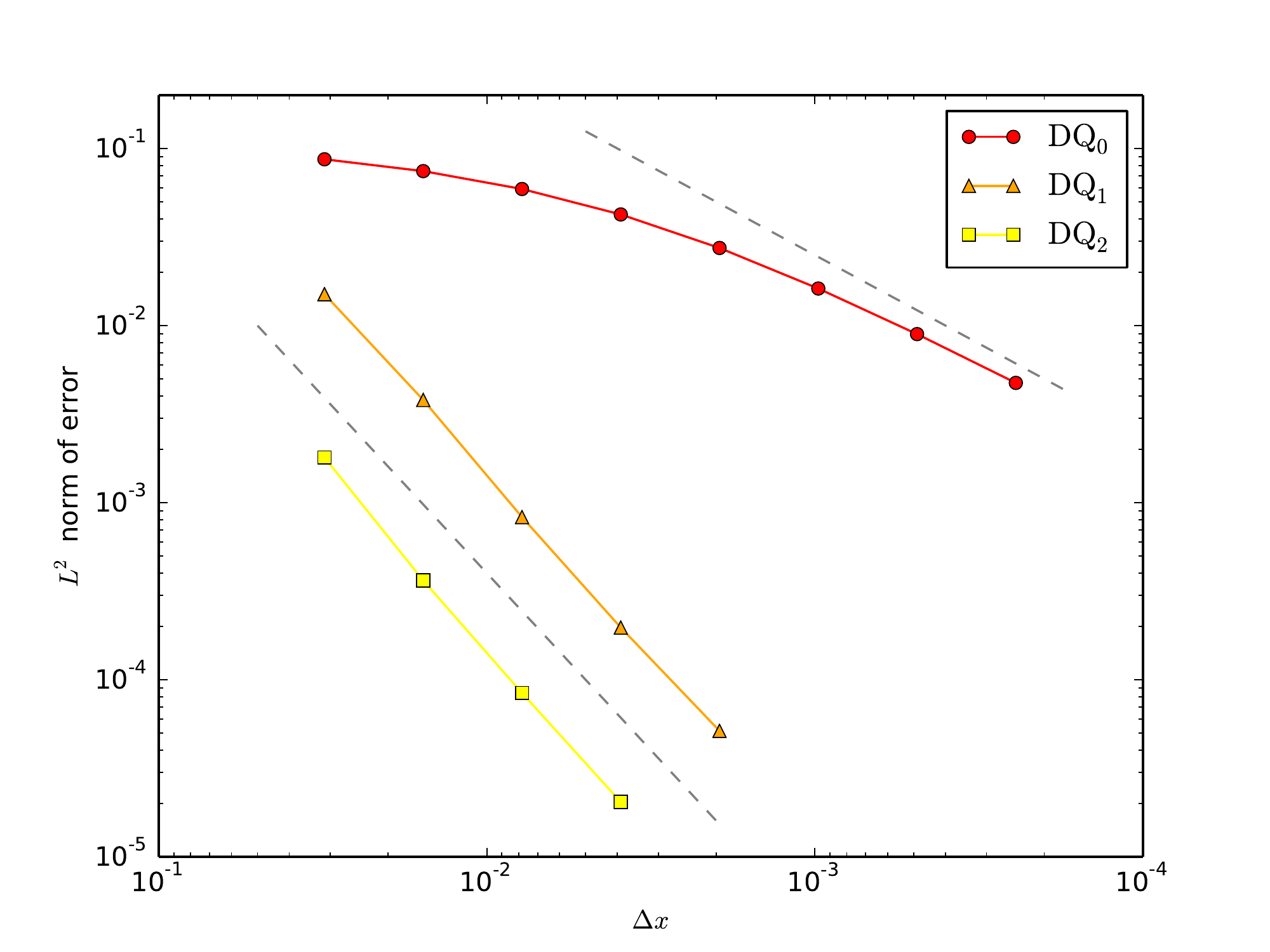}
\caption{The $L^2$ error between the computed and `analytic' solution
is plotted against~$\Delta x$ for the problem described
in \autoref{sec:dgadv-2d}. The dotted lines are proportional
to~$\Delta x$ and~$\Delta x^2$, and are merely to aid comprehension.
The~$\DQ{0}$ simulations converge at first-order for sufficiently
small values of~$\Delta x$. The~$\DQ{1}$ simulations converge at
second-order, as expected. The cosine bell initial condition has a
discontinuous second derivative, which inhibits the~$\DQ{2}$
simulations from exceeding a second-order rate of convergence.}
\label{fig:ex3_plot2d}
\end{figure}

To discretise this problem, we subdivide $\Omega$ into squares with side
length $\Delta x$. We use~$\DQ{r}$ for the discontinuous function space,
for~$r$ from 0 to 2, which are products of 1D discontinuous
elements. We initialise~$q$ by interpolating
the expression given in Eq.~\eqref{eq:dgadv-2d-q-ic} into the appropriate
space. We approximate~$\vec{u_0}$ by interpolating the expression given
in Eq.~\eqref{eq:dgadv-2d-u0} onto a vector-valued function
in~$\Q{2}$. The~$L^2$ errors between the initial and final~$q$
fields for varying $\Delta x$ are plotted in \autoref{fig:ex3_plot2d}.

\section{Limitations and extensions}
\label{sec:limext}

There are several limitations of the current implementation, which
leaves scope for future work. The most obvious is that the quadrature
calculations are relatively inefficient, particularly at high order.
The product structure of the basis functions can be exploited to
generate a more efficient implementation of numerical quadrature. This
can be done using the \emph{sum-factorisation} method, which lifts
invariant terms out of the innermost loop. In the very simplest cases,
direct factorisation of the integral may be possible. Such operations
could have been implemented within Firedrake's form compiler.
However, this would mask the underlying issue -- that FIAT, which is
supposed to be wholly responsible for producing the finite elements, has
no way to communicate any underlying basis function structure. Work is
underway on a more sophisticated layer of software that returns an
\emph{algorithm} for performing a given operation on a finite element,
rather than merely an array of tabulated basis functions.

Firedrake has recently gained full support for non-affine coordinate
transformations. In the \emph{previous} version of the form compiler,
the Jacobian of the coordinate mapping was assumed to be constant
across each cell. This is satisfactory for simplices, since the physical
and reference cells can always be linked by an affine transformation.
However, this statement does not hold for quadrilateral, triangular
prism, or hexahedral cells. Firedrake now evaluates the Jacobian at
quadrature points. This functionality is also necessary for accurate
calculations on curvilinear cells, in which the coordinate
transformation is quadratic or higher-order. This allows, for example,
more faithful representations of a sphere or spherical shell, extending
the work done in \citet{rognes2013automating}.

\section{Conclusion}
\label{sec:conc}

\begin{table}
\renewcommand{\arraystretch}{1.25}
\footnotesize
\centering
\caption{Examples of the construction of standard finite element spaces. In the
left-hand column, we use the notation of the \emph{Periodic Table of the Finite Elements~\citep{arnold2014periodic}} where possible. \label{tbl:elts}}
\begin{tabu}{|p{0.32\linewidth}|c|c|}
\hline
Element & Cell & Construction$^\star$\\
\hline
$\Q{r}$ (also written $\Q{r,r}$) & quadrilateral & $\P{r} \otimes \P{r}$ \\
$\RTCE{r}$, Raviart--Thomas `edge' element$^\dagger$ & quadrilateral & $\texttt{HCurl}(\P{r} \otimes \DP{r-1}) \oplus \texttt{HCurl}(\DP{r-1} \otimes \P{r})$ \\
Nédélec `edge' element of the second kind$^\ddagger$ & quadrilateral & $\texttt{HCurl}(\P{r} \otimes \DP{r}) \oplus \texttt{HCurl}(\DP{r} \otimes \P{r})$ \\
$\RTCF{r}$, Raviart--Thomas `face' element~\citep{raviart1977mixed} & quadrilateral & $\texttt{HDiv}(\P{r} \otimes \DP{r-1}) \oplus \texttt{HDiv}(\DP{r-1} \otimes \P{r})$ \\
Nédélec `face' element of the second kind$^\ddagger$ & quadrilateral & $\texttt{HDiv}(\P{r} \otimes \DP{r}) \oplus \texttt{HDiv}(\DP{r} \otimes \P{r})$ \\
$\DQ{r}$ (discontinuous $\Q{r}$) & quadrilateral & $\DP{r} \otimes \DP{r}$ \\
\hline
$\P{r,r}$$^{\dagger\dagger}$ & triangular prism & $\P{r}^\triangle \otimes \P{r}$ \\
Nédélec `edge' element of the first kind$^{\ddagger\ddagger}$ & triangular prism & $\texttt{HCurl}(\P{r}^\triangle \otimes \DP{r-1}) \oplus \texttt{HCurl}(\RTE{r}^\triangle \otimes \P{r})$ \\
Nédélec `edge' element of the second kind~\citep{nedelec1986new} & triangular prism & $\texttt{HCurl}(\P{r}^\triangle \otimes \DP{r}) \oplus \texttt{HCurl}(\BDME{r}^\triangle \otimes \P{r})$ \\
Nédélec `face' element of the first kind$^{\ddagger\ddagger}$ & triangular prism & $\texttt{HDiv}(\RTF{r}^\triangle \otimes \DP{r-1}) \oplus \texttt{HDiv}(\DP{r-1}^\triangle \otimes \P{r})$ \\
Nédélec `face' element of the second kind~\citep{nedelec1986new} & triangular prism & $\texttt{HDiv}(\BDMF{r}^\triangle \otimes \DP{r}) \oplus \texttt{HDiv}(\DP{r}^\triangle \otimes \P{r})$ \\
$\DP{r,r}$ & triangular prism & $\DP{r}^\triangle \otimes \DP{r}$ \\
\hline
$\Q{r}$ (also written $\Q{r,r,r}$) & hexahedra & $\Q{r}^\square \otimes \P{r}$ \\
$\NCE{r}$, Nédélec `edge' element of the first kind~\citep{nedelec1980mixed} & hexahedra & $\texttt{HCurl}(\Q{r}^\square \otimes \DP{r-1}) \oplus \texttt{HCurl}(\RTCE{r}^\square \otimes \P{r})$ \\
Nédélec `edge' element of the second kind~\citep{nedelec1986new} & hexahedra & $\texttt{HCurl}(\Q{r}^\square \otimes \DP{r}) \oplus \texttt{HCurl}(\mathrm{N2CE}_r^\square \otimes \P{r})$ \\
$\NCF{r}$, Nédélec `face' element of the first kind~\citep{nedelec1980mixed} & hexahedra & $\texttt{HDiv}(\RTCF{r}^\square \otimes \DP{r-1}) \oplus \texttt{HDiv}(\DQ{r-1}^\square \otimes \P{r})$ \\
Nédélec `face' element of the second kind~\citep{nedelec1986new} & hexahedra & $\texttt{HDiv}(\mathrm{N2CF}_r^\square \otimes \DP{r}) \oplus \texttt{HDiv}(\DQ{r}^\square \otimes \P{r})$ \\
$\DQ{r}$ & hexahedra & $\DQ{r}^\square \otimes \DP{r}$ \\
\hline
\end{tabu}

$^\dagger$: this is a curl-conforming analogue of the usual Raviart--Thomas quadrilateral element~\citep{raviart1977mixed}.\\
$^\ddagger$: these are the quadrilateral reductions of the hexahedral Nédélec elements of the second kind~\citep{nedelec1986new}.\\
$^{\dagger\dagger}$: this denotes the element with polynomial degree $r$ in the first two variables, and polynomial degree $r$ in the third variable separately.\\
$^{\ddagger\ddagger}$: these are the prism equivalents of the tetrahedral and hexahedral Nédélec elements~\citep{nedelec1980mixed}.\\
$^\star$: $\RTE{}$ and $\RTF{}$ refer to the Raviart--Thomas edge and face elements on triangles. $\BDME{}$ and $\BDMF{}$ refer to the Brezzi--Douglas--Marini~\citep{brezzi1985two} edge and face elements on triangles. $\mathrm{N2CE}$ and $\mathrm{N2CF}$ refer to the Nédélec elements of the second kind that we construct on quadrilaterals.
\end{table}

This paper presented extensions to the automated code generation
pipeline of Firedrake to facilitate the use of finite element spaces on
non-simplex cells, in two and three dimensions. A wide range of finite
elements can be constructed, including, but not limited to, those listed
in \autoref{tbl:elts}. The examples made extensive use of the recently-added
extruded mesh functionality in Firedrake; a related paper detailing the
implementation of extruded meshes is in preparation.

All numerical experiments given in this paper were performed with the
following versions of software, which we have archived on Zenodo:
Firedrake~\citep{zenodo_firedrake}, PyOP2~\citep{zenodo_pyop2},
TSFC~\citep{zenodo_tsfc}, COFFEE~\citep{zenodo_coffee},
UFL~\citep{zenodo_ufl}, FIAT~\citep{zenodo_fiat},
PETSc~\citep{zenodo_petsc}, PETSc4py~\citep{zenodo_petsc4py}. The
code for the numerical experiments can be found in the supplement
to the paper.

\section*{Acknowledgements}
Andrew McRae wishes to acknowledge funding and other support from the
Grantham Institute and Climate-KIC. This work was supported by the
Natural Environment Research Council [grant numbers NE/K006789/1,
NE/K008951/1], and an Engineering and Physical Sciences Research Council
prize studentship.

\bibliography{paper}

\end{document}